\documentclass[3p]{elsarticle}

\usepackage{lineno,hyperref}
\usepackage{amsmath}    
\usepackage{amssymb,amsthm}    
\usepackage{bm}
\usepackage{wrapfig}
\usepackage{booktabs}
\usepackage{subcaption}
\usepackage{threeparttable}
\usepackage{multirow}
\usepackage{pifont}
\usepackage{units}
\usepackage{float}
\usepackage{blindtext}
\usepackage{soul}
\usepackage{varwidth}
\usepackage{wasysym}
\usepackage{algpseudocode}  
\usepackage{algorithm}      
\algnewcommand\algorithmicto{\textbf{to}}
\algnewcommand\algorithmicin{\textbf{in}}
\algnewcommand\algorithmicforeach{\textbf{for each}}
\algrenewtext{For}[3]{\algorithmicfor\ #1 $\gets$ #2\ \algorithmicto\ #3\ \algorithmicdo}
\algdef{S}[FOR]{ForEach}[2]{\algorithmicforeach\ #1\ \algorithmicin\ #2\ \algorithmicdo}
\algdef{SE}[DOWHILE]{Do}{DoWhile}{\algorithmicdo}[1]{\algorithmicwhile\ #1}%
\modulolinenumbers[2]

\journal{Computer Methods in Applied Mechanics and Engineering}

\newcommand{\dpar}[2]{\ensuremath{\frac{\partial #1}{\partial #2}}}

\newcommand{\R}{\mathbb{R}}
\newcommand{\T}{\mathsf{T}}
\renewcommand{\L}{\mathcal{L}}
\renewcommand{\L}{\mathcal{L}}
\renewcommand{\b}{\boldsymbol}

\newcommand{\x}{\mathbf{x}}
\newcommand{\w}{\mathbf{w}}
\newcommand{\svms}{\sigma_\mathrm{VMS}}
\newcommand{\sxx}{\sigma_{xx}}
\newcommand{\syy}{\sigma_{yy}}
\newcommand{\szz}{\sigma_{zz}}
\newcommand{\sxy}{\sigma_{xy}}
\newcommand{\syz}{\sigma_{yz}}
\newcommand{\sxz}{\sigma_{xz}}

\newcommand{\ezz}{\varepsilon_{zz}}

\newcommand{\eyz}{\varepsilon_{yz}}
\newcommand{\exz}{\varepsilon_{xz}}
\newcommand{\Sph}[1]{B_{#1}}
\newcommand{\Nload}{N_\mathrm{load}}
\newcommand{\Nseg}{N_\mathrm{seg}}
\newcommand\norm[1]{\left \| #1 \right \|}

\newcommand*\mjstar{\ding{73}}
\newcommand*\mjstarsmall{\ding{83}}
\newcommand*\mjtri{$\triangle$}
\usepackage{tcolorbox}
\tcbuselibrary{skins,breakable}
\usetikzlibrary{shadings,shadows}

\newenvironment{block}[3]{%
	\tcolorbox[noparskip,breakable,
		colback=#2!1, 
		colframe=#2!20, 
		coltitle=#3, 
		title=\textbf{#1},
	]}%
{\endtcolorbox}

\graphicspath{{./figures/}}

\begin{document}

\begin{frontmatter}

	\title{Numerical analysis of small-strain elasto-plastic deformation using local Radial Basis Function approximation with Picard iteration}

	\author[ijs]{Filip Strniša}
	\ead{filip.strnisa@ijs.si}
	\author[ijs]{Mitja Jančič}
	\ead{mitja.jancic@ijs.si}
	\author[ijs]{Gregor Kosec\corref{correspondingauthor}}
	\ead{gregor.kosec@ijs.si}

	\cortext[correspondingauthor]{Corresponding author}

	\address[ijs]{``Jožef Stefan'' Institute, Parallel and Distributed Systems Laboratory, Jamova cesta 39, 1000 Ljubljana, Slovenia}

	\begin{abstract}
		This paper deals with a numerical analysis of plastic deformation under various conditions, utilizing Radial Basis Function (RBF) approximation. The focus is on the elasto-plastic von Mises problem under plane-strain assumption. Elastic deformation is modelled using the Navier-Cauchy equation. In regions where the von Mises stress surpasses the yield stress, corrections are applied locally through a return mapping algorithm. The non-linear deformation problem in the plastic domain is solved using the Picard iteration.
		The solutions for the Navier-Cauchy equation are computed using the Radial Basis Function-Generated Finite Differences (RBF-FD) meshless method using only scattered nodes in a strong form. Verification of the method is performed through the analysis of an internally pressurized thick-walled cylinder subjected to varying loading conditions. These conditions induce states of elastic expansion, perfectly-plastic yielding, and plastic yielding with linear hardening. The results are benchmarked against analytical solutions and traditional Finite Element Method (FEM) solutions. The paper also showcases the robustness of this approach by solving case of thick-walled cylinder with cut-outs. The results affirm that the RBF-FD method produces results comparable to those obtained through FEM, while offering substantial benefits in managing complex geometries without the necessity for conventional meshing, along with other benefits of meshless methods.
	\end{abstract}

	\begin{highlights}
		\item The elasto-plastic von Mises problem is solved in strong form on scattered nodes using direct RBF-FD and Picard iteration.
		
		\item The presented method provides results regardless the underlying geometry without any special treatment of differential operators, boundary conditions or spatial discretisation.
		
		\item The presented method shows good agreement with analytical and Finite Element Method solutions, as well as convergence behaviour.

		\item The method effectively handles different stages of material behavior, from elastic deformation through perfect plastic yielding to linear hardening conditions.

	\end{highlights}


\end{frontmatter}


\section{Introduction}
The study of plasto-elastic behaviour is a field of research concerned with the challenges associated with modelling the deformation of materials that exhibit both elastic and plastic behaviour under applied loads~\cite{chakrabarty1987, hill1998}. These materials, which are widely used in engineering and manufacturing, initially exhibit linear elastic behaviour when subjected to mechanical stress, but begin to deform plastically as soon as a threshold yield stress is exceeded. Accurate prediction of this transition is important for the design and analysis of structural components subjected to dynamic loads such as earthquakes, for modelling the behaviour of soils and rocks under construction loads or natural forces, and for additive manufacturing processes where materials are often subjected to complex loads, to name but a few applications.

Elastic deformations are usually modelled using the Navier-Cauchy equation. As soon as the material begins yielding (i.e. the stress locally exceeds the yield stress), the Navier-Cauchy solution is considered as an initial guess in an iterative algorithm that adjusts the solution locally.
Traditionally, the finite element method (FEM) is used to solve elasto-plastic problems~\cite{desouza2008,Schrder2015SmallSP, Roostaei2018ACS, amouzou2021}. Although FEM is a powerful tool that offers a mature and versatile solution approach including all types of adaptivities~\cite{mitchell2014comparison}, well-understood error indicators~\cite{segeth2010review} and an established toolset for isogeometric analysis, i.e. coupling with computer-aided design (CAD)~\cite{cottrell2009isogeometric}, the complexity of mesh generation has led researchers to explore alternative methods. In FEM, the domain is discretised by dividing a complex geometry into smaller, simpler parts, called elements, that cover it enitely. Despite significant developments in the field of mesh generation, the process of meshing often remains the most time-consuming part of the overall solution procedure, while the quality of the meshing has a direct impact on the accuracy, convergence and computational cost of the method~\cite{liu_introduction_2005}.

In contrast, meshless methods work with node clouds, which are often referred to as scattered nodes and do not require a topological relationship between them. Although some authors argued that arbitrary nodes could be used~\cite{liu2002mesh, Reuther2012}, nowadays it is generally accepted that for a stable meshless approximation the domain should be discretised with quasi-uniformly~\cite{wendland2004scattered,liu2009meshfree} distributed scattered nodes. Modern node positioning algorithms can automatically populate complex geometries, are dimension-independent, support variable node densities and can handle domains whose boundaries are defined by CAD objects~\cite{slak2019generation,duh2021fast,fornberg2015fast,van2021fast,shankar2018robust, duh2023discretization}, and yet are significantly less complex than meshing~\cite{slak2019generation}.

The meshless methods have already been applied in elasto-plastic problems. In~\cite{ji2005meshfree} authors demonstrated the Reproducing Kernel Particle Method weak solution of Drucker-Prager and Mohr-Coulomb models. Another weak form solution was presented in~\cite{kargarnovin2004elasto}, where authors demonstrated the Element free Galerkin method in elasto-plastic stress analysis around the tip of a crack. In~\cite{app132312591} a Radial Point Interpolation method was demonstrated in elasto-plastic analysis of frame structures.

In weak form methods, the governing equations of plasticity are formulated in an integral form rather than a differential form. The key difference from strong form meshless methods, which directly discretize the differential equations at nodes, is that weak form methods satisfy these equations in an average sense over the domain, hence potentially increasing the solution's smoothness and numerical stability. However, a recent study revealed that strong form solution might perform better in capturing the peaks in the stress comparing it to the weak from solutions~\cite{kosec_weak_2019}.

In a strong form, elasto-plastic deformation was approached with the meshless method of fundamental solutions (MFS) in the simulation of material behaviour that hardens with plastic deformation characterized by the Chakrabarty model~\cite{JANKOWSKA201812}. In the RBF-FD context, the problem was recently addressed in the analysis of cantilever beam and the Reinforced Concrete column-Steel~\cite{jiang2021nonlinear}. In~\cite{jiang2021nonlinear} authors used a multiquadratic basis that requires shape tuning, which can have a notable effect on the accuracy of the numerical solution due to stagnation errors and on the stability of the solution procedure~\cite{bayona2017role, wang2002optimal}. Solution presented in~\cite{jiang2021nonlinear} also stabilises the calculation of the first derivatives using Finite difference method (FDM) on bilinear local nodes. Since the study is limited to regular nodes, authors only discussed how to extend the approach to scattered nodes by using fictitious nodes at which the field nodes should be interpolated. This has been properly implemented in the recent paper~\cite{vuga_improved_2024}, where authors argued that direct RBF-FD, i.e. without using FDM to compute the divergence operator, is unsuitable for solving elasto-plastic problems. Furthermore, authors introduced a stabilisation of the boundary conditions by adding the first two layers of nodes inside the domain on the opposite side of the outward-facing normal vector instead of using purely scattered nodes. They have shown that such stabilisation is essential for the convergence of a method. Although these stabilisations improve the stability of the method, they also reduce to some extent the elegance of the meshless discretisation of the domain and introduce new free parameters that need to be tuned.

In this paper, we present an alternative approach that does not require any kind of stabilisation, i.e. a direct RBF-FD can be applied to purely scattered nodes. In~\cite{vuga_improved_2024}, the authors solve the system of nonlinear equations arising from plastic deformation using the Newton–Raphson iteration. Here, we propose to solve the problem using the Picard iteration, a much simpler approach where there is no need to assembly the Jacobian matrix and only elastic material parameters are used, resulting in a constant stiffness matrix during the iteration. Although the Picard iteration converges significantly slower than the Newton–Raphson iteration~\cite{yarushina2010}, it is less computationally expensive per iteration and, as shown in this paper, 
converges without stabilisation of the divergence operator or boundary conditions. The application of Picard iteration in field of elasto-plastic materials was initially proposed by Kołodziej et al.~\cite{KOLODZIEJ20134217}, where authors considered elasto–plastic properties of prismatic bars using the MFS. The use of Picard iteration was also recently reported in elasto-plastic torsion problem~\cite{moayyedian2021elastic} using the MFS as well as the Taylor series expansion with moving least squares in~\cite{XU2023939}.

In this paper, we discuss a von Mises plasticity model with non-linear isotropic hardening under small-strains assumption in a plane stress example of internally pressurized thick-walled cylinder subjected to varying loading conditions. The solution procedure is built on the direct RBF-FD approximation with Polyharmonic Splines (PHS) basis functions and polynomial augmentation and Picard iteration. The implementation was done using our in-house Medusa C++ library~\cite{slak2021medusa} supporting all the required meshless procedures. 

\section{Iterative algorithm for solving Elasto-Plastic deformation problems}
\label{sec:plasticity}
\begin{figure} [H]
	\centering
	\includegraphics[width=0.5\columnwidth]{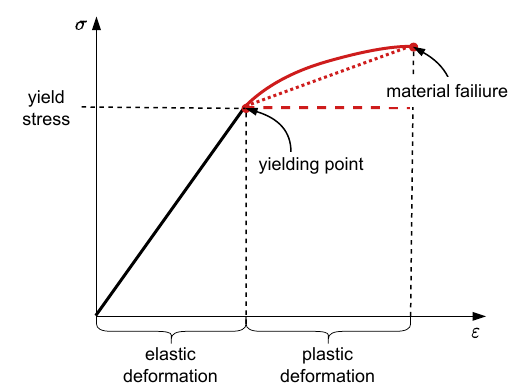}
	\caption{Sketch of stress and strain relation as a result of an uniaxial tension test.}
	\label{fig:stress_strain_sketch}
\end{figure}

Plastic deformation, also known as \emph{plasticity}, is the phenomenon when a shape of a solid material undergo non-reversible changes as a response to the applied forces, i.e.\ material's internal structure is changed during the deformation. The transition from elastic to plastic behaviour is known as \emph{yielding} and the ability to predict it is important for the engineering of structures and components that are subjected to various loads and exposed to diverse environmental conditions.
Typically, uniaxial tension tests are conducted to gain insights into the response characteristics of a material when exposed to external loads. Such tests, when performed with a ductile material, yield a \emph{stress-strain} relationship curve as conceptually illustrated in Figure~\ref{fig:stress_strain_sketch}.

In the first segment, the material undergoes an \emph{elastic deformation}. Since the stress-strain relation in that segment is linear, it is often referred to as \emph{linear elastic deformation}, where the slope of the curve is the elastic Young's modulus $E$. If the applied load is released before the \emph{yield stress} is achieved, the material returns to its original state without permanent structural changes. On the other hand, if applied loads introduce internal stresses exceeding the yield stress, the total observed deformation becomes a sum of two contributions --- the plastic part and the elastic part.

The stress-strain relation beyond the \emph{yielding point} is commonly referred to as the \emph{hardening curve}. When the material under consideration is unloaded in that segment, the object no longer returns to its original shape. Instead, at least parts of it, stay \emph{plastically deformed}. Different materials exhibits different hardening curves, marked with red in Figure~\ref{fig:stress_strain_sketch}. A commonly considered simplification of the hardening curve is \emph{ideal plasticity}, where the yield stress beyond the yielding point is constant (see dashed red line in Figure~\ref{fig:stress_strain_sketch}). Another common modelling approach is to use linear hardening curve marked with red dotted line. Such modelling approach results in a \emph{bilinear} elastic-plastic material behaviour, because both elastic and plastic regimes are governed by a linear stress-strain relation with the Young's modulus $E$ denoting the slope in the elastic regime and hardening rate $H$ denoting the slope in the plastic regime. In practice, the hardening is commonly represented as curve fitted to experimental data obtained from uniaxial tension tests, in simplest case a piecewise linear fit between the experimental data resulting in \emph{multilinear material models}.

The study presented in this work is limited to rate-independent materials where the deformation of the material under consideration does not depend on the rate of the applied loads. The research is also limited to \emph{isotropic hardening}, where the evolution of yield stress surface corresponds to a uniform (isotropic) expansion of the initial yield surface. A detailed discussion on other stress-strain relationships and different types of hardening can be found in~\cite{desouza2008}.
%

Initially, the deformation is elastic, as depicted in Figure~\ref{fig:stress_strain_sketch}.
The principal assumption of linear elasticity is that, when subjected to external forces, a body deforms linearly -- obeying Hooke's law. Thus, the stress tensor $\bm{\sigma}$ is defined as the product of the elasticity tensor $\mathrm{D}^\mathrm{e}$ and the strain tensor $\bm{\varepsilon}$,
\begin{equation}
	\label{eq: stress}
	\bm{\sigma} = \mathrm{D}^\mathrm{e} \bm{\varepsilon}.
\end{equation}
commonly expressed as
$\mathrm{D}^\mathrm{e} = 2\mu \mathrm{I}_S + \lambda\mathbf{I} \otimes \mathbf{I}$,
where $\mathrm{I}_S$ is the symmetric identity tensor, $\mathbf{I}$ is the identity matrix and $\lambda$ and $\mu$ are the first and second Lamé parameters, respectively, expressed in terms of Young's modulus $E$, and Poisson's ratio $\nu$ as
\begin{equation}
	\label{eq: lame}
	\lambda = \frac{E \nu}{\left(1 - 2\nu\right)\left(1 + \nu\right)} \quad \text{and} \quad \mu     = \frac{e}{2 \left(1 + \nu\right)}.
\end{equation}
Under small-strain assumption the strain tensor $\bm \varepsilon$ and displacement vector $\vec{u}$ are related as
\begin{equation}
	\label{eq: strain}
	\bm{\varepsilon} = \frac{\nabla \vec{u} + (\nabla \vec{u})^{\intercal}}{2}.
\end{equation}
Considering \eqref{eq: stress} and \eqref{eq: strain} and the stationary Cauchy momentum equation in the absence of the external force ($\nabla \cdot \bm{\sigma}  = 0$), one gets a well-known Navier-Cauchy equation
\begin{equation}
	\label{eq: navier-cauchy}
	\left(\lambda + \mu\right) \nabla (\nabla \cdot \vec{u}) + \mu \nabla^2 \vec{u} = 0.
\end{equation}
describing the linear-elastic deformation.
Two main types of boundary conditions are frequently employed: essential boundary conditions and traction boundary conditions, also referred to as natural boundary conditions. Essential boundary conditions are used to define displacements $\vec{u}_0$ along certain parts of the domain's boundary, denoted as $\vec{u} = \vec{u}_0$. On the other hand, traction boundary conditions are concerned with defining surface traction $\vec{t}_0$, expressed as $\bm{\sigma}\vec{n} = \vec{t}_0$, were, $\vec{n}$ represents the outward unit normal to the domain's boundary. 

Upon surpassing a critical stress threshold (\emph{yield stress}~$\sigma_y$), the material experiences localized irreversible plastic deformation. In the plastic regime, the total local deformation has two contributions. One due to the elastic strain $\bm{\varepsilon}^\mathrm{e}$ and another due to the plastic strain $\bm{\varepsilon}^p$. The sum of both is referred to as the total strain $\bm{\varepsilon} = \bm{\varepsilon}^e + \bm{\varepsilon}^p$. In plastic regime the Navier-Cauchy~\eqref{eq: navier-cauchy} is therefore not valid anymore. To check its validity different criteria as the Tresca, the Mohr-Coulomb, the von Mises yield or the Drucker-Prager criteria~\cite{desouza2008} can be used. In this paper we focus on the von Mises yield criterion assuming that the yield stress ($\sigma_y$) is a function of the scalar equivalent plastic strain $\overline\varepsilon^p$ (see Figure~\ref{fig:irregular_stress_strain} for example of such a stress-strain relationship). The criterion is based on the the yield function $\Phi$ defined as
\begin{equation}
	\label{eq: y funct}
	\Phi = {\svms}_n - \sigma_y(\overline\varepsilon^p_n),
\end{equation}
where von Mises stress $\svms$ stands for
\begin{equation}
	\label{eq: vms}
	\begin{split}
		\svms &= \frac{1}{\sqrt{2}}\left[\left(\sxx - \syy\right)^2 + \left(\syy - \szz\right)^2 + \left(\szz - \sxx\right)^2\right. \\
		&\left.\phantom{\left(\sxx - \syy\right)^2}\quad + 6 \left(\sxy^2 + \syz^2 + \sxz^2\right)\right]^\frac{1}{2}. 
	\end{split}
\end{equation}
If $\Phi > 0$, the material is locally under plastic deformation, otherwise the deformation is elastic
\begin{equation}
	\label{eq: criterion}
	\left\{\begin{matrix}
		\Phi \leq 0, & \quad \text{(elastic deformation), solution is locally valid}   \\
		\Phi > 0   , & \quad \text{(plastic deformation), solution is locally invalid.}
	\end{matrix}\right. 
\end{equation}
In the elastic regime, the Navier-Cauchy solution is accepted, while in the plastic regime a correction is applied that adjusts the calculated stresses and strains to account for plastic effects.

Since the state of stress reflects the history of applied stress, the problem can be solved in incremental way by applying partial loads to the body 

that eventually add-up to the total applied load~\cite{desouza2008,yarushina2010}. The key concept behind such approach is that every partial load applied to the body is assumed to result in a purely elastic deformation. Thus, a first estimate is obtained by solving the Navier-Cauchy Equation~\eqref{eq: navier-cauchy} resulting in a vector field of displacement increments $\delta\bm{u}$ corresponding to the amount of applied partial load. In the $i$-th load step, $i\leq \Nload$, one obtains the following set of bulk state predictions
\begin{equation}
	\label{eq: pred}
	\begin{aligned}
		\vec{u}_{i}                     & = \vec{u}_{i - 1} + \delta\vec{u}                                                               \\
		\delta\bm{\varepsilon}          & = \frac{\nabla \left(\delta\vec{u}\right) + (\nabla \left(\delta\vec{u}\right))^{\intercal}}{2} \\
		\bm{\varepsilon}^\mathrm{e}_{i} & = \bm{\varepsilon}^\mathrm{e}_{i - 1} + \delta\bm{\varepsilon}                                  \\
		\bm{\varepsilon}^p_{i}          & = \bm{\varepsilon}^p_{i - 1}                                                                    \\
		\overline \varepsilon^p_i       & = \overline \varepsilon^p_{i - 1}                                                               \\
		\bm{\sigma}_i                   & = \mathrm{D}^\mathrm{e} \bm{\varepsilon}^\mathrm{e}_i
	\end{aligned}
\end{equation}
where a scalar value $\overline \varepsilon^p$ denotes the accumulated plastic strain, and  the stress tensor is calculated only from the elastic strain. In next step, we check the validity of the stress by assessing the yield criterion~\eqref{eq: criterion}.
If $\Phi \leq 0$, the solution is valid as the yield condition has not been violated, and the assumed elastic deformation computed using~\eqref{eq: navier-cauchy} holds. If $\Phi > 0$, on the other hand, the prediction was incorrect, thus, solution needs to first be locally adjusted according to the \emph{return mapping procedure}, where we use the single-equation return mapping scheme~\cite{desouza2008}. For this purpose an incremental plastic multiplier $\Delta\gamma$ is introduced. The goal of the return mapping procedure is to find such $\Delta\gamma$, that the following equality is satisfied~\cite{desouza2008}
\begin{equation}
	\label{eq: rm goal}
	\Phi = {\svms}_i  - 3 \mu \Delta\gamma - \sigma_y(\overline\varepsilon^p_i + \Delta\gamma) = 0.
\end{equation}
We use a Newton-Raphson approach (see Algorithm~\ref{alg:return_mapping} for details) to solve~\eqref{eq: rm goal} and finally, with the computed $\Delta\gamma$, the local bulk state variables are temporarily updated
\begin{equation}
	\label{eq: updated}
	\begin{aligned}
		\mathbf s _{i}^{j+1}                  & = \bigg ( 1 - \frac{\Delta \gamma 3 \mu}{{\svms}_i^j}\bigg ) \mathbf s_i^j                 \\
		\bm\sigma_{i}^{j+1}                   & = \mathbf s_{i+1}^{j+1} + p_{i+1}^{j+1}\mathbf I                                              \\
		\bm{\varepsilon}^{\mathrm{e} \ j+1}_{i} & = \frac{1}{2\mu}\mathbf{s}_{i + 1}^{j+1} + \frac{1}{3}\varepsilon^{\mathrm{e} \ j}_{v \ i} \mathbf I \\
		\overline \varepsilon^{p \ j+1}_{i}     & = \overline \varepsilon^{p \ j}_{i} + \Delta \gamma
	\end{aligned}
\end{equation}
where $\mathbf s$ and $p$ denote the deviatoric and hydrostatic components of the stress tensor and $\varepsilon^\mathrm{e}_{v} = \sum\mathrm{diag}(\bm \varepsilon^\mathrm{e})$ is the volumetric component of the elastic strain and $j$ denotes the global optimisation index (explained in following few sentences). The updated stress $\bm\sigma_{i}^{j+1}$ does not fulfill Navier-Cauchy equation anymore, i.e. 
\begin{equation}
 	\nabla \cdot \bm\sigma_{i}^{j+1} = -\delta \vec{r} ,
 	\label{eq: div_free_iter}
\end{equation}
where vector $\vec{r}\neq 0$ stands for residual. To drive the state towards correct solution a global optimisation is performed, where the whole process repeats using~\eqref{eq: div_free_iter} to calculate stress guess until
\begin{equation}
	 \norm{\delta \vec{r}} < \epsilon_{\mathrm{tol}}
\end{equation}
holds in all nodes with $\epsilon_{\mathrm{tol}}$ standing for prescribed tolerance~\cite{yarushina2010}. Note, that during global optimization, the boundary conditions need to be updated accordingly, e.g., at the $i$-th load step, the total prescribed stress at the boundary $\Omega$ is 
\begin{equation}
	\bm \sigma(\Omega) = \bm \sigma_{\mathrm{load} \ i} = \sum_i \delta \bm \sigma,
\end{equation}
and after $j$ global optimization steps, this boundary condition reads 
\begin{equation}
	\bm \sigma(\Omega) = \bm \sigma_{\mathrm{load} \ i} - \bm \sigma(\Omega)_i^j.
\end{equation}

Once the current load step converges, a new partial load is applied, and the whole process is repeated until the sum of partial loads equals the prescribed total load. 

\begin{algorithm}[h]
	\caption{Return mapping algorithm.}
	\label{alg:return_mapping}
	\vspace{1mm}
	\textbf{Input:} Material properties $\Pi$, maximum return mapping iterations $I_\mathrm{max}$ and return mapping tolerance $\epsilon_\mathrm{tol}$.\\
	\textbf{Output:} Incremental plastic multiplier $\Delta \gamma$.
	\begin{algorithmic}[1]
		\Function{return\_mapping}{$\Pi, I_\mathrm{max}, \epsilon_\mathrm{tol}$}
		\State $k \gets 0$
		\Comment{Initializing iteration counter.}
		\State $\Phi^\star \gets {\svms}_n - \sigma_y(\overline\varepsilon^p_n)$
		\Comment{Initializing residual yield function value.}
		\newline
		\State // Newton-Raphson iteration:
		\Do
		\State $H \gets \frac{\mathrm{d} \sigma_y}{\mathrm{d}\overline{\varepsilon}^p} \Big \rvert_{\overline{\varepsilon}^p_n + \Delta \gamma}$
		\Comment{Hardening curve slope evaluated at $\overline{\varepsilon}^p_n + \Delta \gamma$.}
		\State $\Delta\gamma \gets \Delta \gamma + \frac{\Phi^\star}{3\mu + H}$
		\Comment{New guess for $\Delta\gamma$.}
		\State $\Phi^\star = {\svms}_n  - 3 \mu \Delta\gamma - \sigma_y(\overline\varepsilon^p_n + \Delta\gamma)$
		\Comment{Re-evaluate residual.}
		\State $k = k+ 1$
		\Comment{Increase number of Newton-Raphson iterations.}
		\DoWhile{$\lvert \Phi^\star \rvert > \epsilon_{\mathrm{tol}} \text{ and } k < k_\mathrm{max}$}
		\newline
		\State \Return $\Delta \gamma$
		\Comment{Return incremental plastic multiplier $\Delta \gamma$.}
		\EndFunction
	\end{algorithmic}
\end{algorithm}

The backbone of above solution procedure is the numerical treatment of Navier-Cauchy~\eqref{eq: navier-cauchy} PDE that requires appropriate domain discretisation, partial differential operators discretisation and ultimately solution of global sparse system resulting from the discretisation of the~\eqref{eq: navier-cauchy}. All these steps will be discussed in details in following section. 

The entire solution procedure demonstrated in Algorithm~\ref{alg:adapt}.

\begin{algorithm}[h]
	\caption{Solution procedure algorithm.}
	\label{alg:adapt}
	\vspace{1mm}
	\textbf{Input:} The problem, nodal density function $h$, stencil size $n$, linear differential operators $\L$, approximation basis $\xi$, material properties $\Pi$, number of load steps $N_{\mathrm{load}}$,  external load $\sigma_{\mathrm{load}}$, maximum return mapping iterations $I_\mathrm{max}$ and return mapping tolerance $\epsilon_\mathrm{tol}$.\\
	\textbf{Output:} Deformation field ${\mathbf u}$ and corresponding stress and strain tensors $\b \sigma$ and $\b \varepsilon$ respectively.
	\begin{algorithmic}[1]
		\Function{solve}{$\mathrm{problem}, h, n, \xi, \Pi, N_{\mathrm{load}}, \sigma_{\mathrm{load}}, I_\mathrm{max}, \epsilon_\mathrm{tol}$}
		\State $\Omega \gets \Call{discretise\_domain}{h}$
		\Comment{Discretize the domain.}
		\State $\L \gets \Call{aproximate}{\mathrm{problem}, \xi, \Omega, n}$
		\Comment{Differential operator approximation.}
		\For{$i$}{0}{$N_{\mathrm{load}}$} \label{alg:while}
		\State $\delta \vec u = \mathbf 0$
		\Comment{Initalize displacement changes.}
		\State $\delta \vec r = \mathbf 0$
		\Comment{Initalize residuum force changes.}
		\\
		\Do
		\State $\delta \vec u \gets \Call{solve\_elastic}{\text{problem}, \Pi, \delta \vec r, \partial \sigma_{\mathrm{load}}}$
		\Comment{\parbox[t]{.35\linewidth}{Linear-elastic guess for partial load $\partial \sigma_{\mathrm{load}} = \frac{\sigma_{\mathrm{load}}}{N_{\mathrm{load}}}$.}}\label{alg:partial_load}
		\State $\vec u = \vec u + \delta \vec u$
		\State $\bm \varepsilon,\bm \sigma \gets \Call{process}{\vec u, \Pi}$
		\Comment{Compute strain and stress tensors.}\\
		\ForEach{$\vec p$}{$\Omega$}
		\If{$\svms(\vec p) > \sigma_{\mathrm{yield}}$}~\label{alg:check}
		\State $\bm \sigma (\vec p),\bm \varepsilon(\vec p) \gets \Call{return\_mapping}{I_{\mathrm{max}}, \epsilon_{\mathrm{tol}}, \Pi}$
		\EndIf
		\EndFor
		\State $\delta \vec r \gets \Call{compute\_residuum}{\bm \sigma}$
		\Comment{Estimate residuum forces.}
		\DoWhile{$\Call{max}{\norm{\delta \vec r}} >\epsilon_{\mathrm{tol}}$}
		\Comment{Repeat until residuum forces are sufficiently small.}
		\EndFor
		\State \Return $\bm \varepsilon,\bm \sigma, \vec u$
		\Comment{Return strains, stresses and displacements for all computational nodes.}
		\EndFunction
	\end{algorithmic}
\end{algorithm}

\section{Solution to Navier-Cauchy equation using RBF-FD approximation}
\label{sec:meshless}
In this section we will discuss how to numerically solve the Navier-Cauchy equation~\eqref{eq: navier-cauchy} using the RBF approximation on scattered nodes. In the first step, the domain $\Omega$ and its boundary are populated with discrete nodes at which the solution will be computed. Afterwards, a linear differential operators governing the equation at hand are approximated at each node $\x_i \in \Omega$. In the third step, the Navier-Cauchy equation is approximated with system of algebraic equations that is represented as a global sparse system. Finally, in the last step, the sparse system is solved, resulting in the numerical solution. In following subsections each of the above steps are discussed.

\subsection{Node generation}
\label{sec:domain_discretization}
In RBF community, a well established approach towards stable and accurate RBF approximation is to discretise domain of interest with quasi-uniformly~\cite{wendland2004scattered,liu2009meshfree} scattered nodes. Although there is still no generally accepted consensus on measure for node quality, there are two quantities that are often considered for that purpose, namely the minimum spacing between any pair of nodes that is often also referred to as separation distance, and the maximal empty space without nodes commonly also referred to as fill distance. Here we use an algorithm that generates nodes according the above guidelines~\cite{slak2019generation} and has been used several times in RBF related research~\cite{janvcivc2021monomial, janvcivc2023strong, najafi_divergence-free_2022, berljavac_rbf-fd_2021, jancic_meshless_2024}.¸

The algorithm is of an iterative advancing-front type. It operates by starting with seed nodes that are put in an \emph{expansion queue}. In our particular case, we initially distribute nodes on the boundary, and since the boundary is one-dimensional, they can trivially be distributed uniformly (left-most plot of Figure~\ref{fig:elast_disp_conv}). In each iteration a node is selected from a queue and expanded by generating new candidate (in our case $12$ candidates) nodes uniformly on a randomly rotated annulus around the processed point. Note that for different randomizations, i.e. different random \emph{seeds}, of the annulus rotations we get different node layouts that all conform to the quasi-uniformity guidelines.
These candidates are then processed, discarding those outside the domain or too close to existing nodes. In last iteration step, the nodes that were not discarded are put in a queue. The iteration continues in a breadth-first-search manner until the queue is empty, ensuring the domain is fully covered with the nodes. The algorithm is dimension-independent and capable of populating complex geometries. It also allows for variable internodal distances that can be used as a cornerstone for $h$~\cite{slak2019rbffd} or $hp$~\cite{janvcivc2023strong} adaptivities. It has been also promoted to generation of nodes on parametric surfaces~\cite{duh2021fast} as well as domains whose boundaries are represented as CAD models~\cite{duh2023discretization}, and has been tailored for execution in parallel on shared memory architectures~\cite{depolli_parallel_2022}. Interested readers are referred to the original paper~\cite{slak2019generation} for more details on the node generation algorithm or its stand-alone C++ implementation in the \emph{Medusa library}~\cite{slak2021medusa}.

\subsection{Approximation of Differential Operators}
\label{sec:rbffd}
After computational nodes have been positioned, linear differential operators of the governing system of PDEs are approximated using a set of nearby nodes, commonly referred to as \emph{stencil nodes}. There have been different stencil strategies proposed~\cite{davydov2011adaptive,davydov2023improved}, nevertheless, we resort to the most common approach using cloud of closest $n$ nodes. 
Let us assume a linear differential operator $\L$ in point $\x _c \in \Omega$ and its stencil nodes $\left \{ \x_c \right\}_{i=1}^n = \mathcal{N}$. The approximation is then sough using ansatz
\begin{equation}
	\label{eq:ansatz}
	(\L g)(\vec{x}_c) \approx \sum_{i=1}^n w^{\L}_ig(\vec{x}_i),
\end{equation}
for any function $g$ and unknown weights $\mathbf{w}^{\L}$ that are obtained by enforcing the equality of approximation~\eqref{eq:ansatz} for a chosen set of basis functions $\mathbf{f}(\vec{x})$.

\begin{equation}
	\sum_{i=1}^nw_if_j(\vec{x_i})=(\L f_j)(\vec{x_c}).
\end{equation}
Gathering all $n$ equations for a given stencil results in a system
\begin{equation}
	\label{eq:system_rbf}
	\mathbf F\w =\b{\ell}_f, \qquad
	\mathbf F = \begin{bmatrix}
		f(\norm{\vec{x_1} - \vec{x_1}}) & \cdots & f(\norm{\vec{x_n} - \vec{x_1}}) \\
		\vdots                          & \ddots & \vdots                          \\
		f(\norm{\vec{x_1} - \vec{x_n}}) & \cdots & f(\norm{\vec{x_n} - \vec{x_n}}) \\
	\end{bmatrix}, \qquad
	\ell_f^i = \Big(\L f_i(\vec{x})\Big)\Big|_{ \vec{x} = \vec{x_c}}.
\end{equation}

In this work, Polyharmonic Splines (PHS)  are used as an approximation basis defined as
\begin{equation}
	f(r) = \begin{cases}r^k,       & k \text{ odd}  \\
             r^k\log r, & k \text{ even}\end{cases}
\end{equation}
for Euclidean distance $r$, the dependency on the shape parameter is eliminated~\cite{bayona2019comparison}.
It has previously been observed that using only PHS as basis functions does not guarantee convergent behaviour or solvability~\cite{wendland2004scattered}. Convergent behaviour and conditional positive definiteness can be assured if the approximation basis is expanded with monomials~\cite{bayona2019insight, wendland2004scattered, flyer2016role}. Therefore, the system~\eqref{eq:system_rbf} is augmented with $N_p = \binom{m+d}{d}$ monomials $p$ with orders up to (and including) degree $m$ for a $d$-dimensional domain, yielding an overdetermined system of equations describing the constrained optimization problem~\cite{flyer2016}. In practice, the weights $\w$ are expressed as a solution to the
\begin{equation} \label{eq:rbf-system-aug}
	\begin{bmatrix}
		\mathbf F    & \mathbf P \\
		\mathbf P^\T & \mathbf 0
	\end{bmatrix}
	\begin{bmatrix}
		\mathbf w \\
		\b \lambda
	\end{bmatrix}
	=
	\begin{bmatrix}
		\b\ell_f \\
		\b\ell_p
	\end{bmatrix},
\end{equation}
with the matrix $\mathbf P$ of the evaluated monomials and Lagrangian multipliers $\mathbf \lambda$ that are discarded once the weights are computed.

The discussed approach is also known as a RBF generated finite differences method (RBF-FD) that is a popular choice in the meshless community due to its stability~\cite{le2023guidelines,9803334}.

\subsection{Navier-Cauchy equation Discretization}
In previous sections we discussed how to numerically approximate differential operators on scattered nodes, in this subsection, we discuss the discretisation of the Navier-Cauchy equation~\eqref{eq: navier-cauchy}. Since we deal with a two dimensional domain, we first express the governing problem component-wise $\vec{u} = (u, v)$ as
\begin{align}
	(\lambda+\mu) \frac{\partial}{\partial x}\left( \frac{\partial u}{\partial x} + \frac{\partial v}{\partial y} \right) + \mu \left( \frac{\partial^2 u}{\partial x^2} + \frac{\partial^2 u}{\partial y^2} \right) & = 0 \label{eq:navier-u}  \\
	(\lambda+\mu) \frac{\partial}{\partial y}\left( \frac{\partial u}{\partial x} + \frac{\partial v}{\partial y} \right) + \mu \left( \frac{\partial^2 v}{\partial x^2} + \frac{\partial^2 v}{\partial y^2} \right) & = 0 \label{eq:navier_2d}
\end{align}
and natural boundary conditions
\begin{align}
	{t_0}_x & = \mu n_y\dpar{u}{y} + \lambda n_x \dpar{v}{y} + (2\mu+\lambda)n_x\dpar{u}{x} + \mu n_y \dpar{v}{x} \label{eq:trac1} \\
	{t_0}_y & = \mu n_x\dpar{u}{y} + (2\mu+\lambda) n_y \dpar{v}{y} + \lambda n_2\dpar{u}{x} + \mu n_x \dpar{v}{x}
	\label{eq:trac2}
\end{align}
with surface traction $\vec{t_0} = ({t_0}_x, {t_0}_x)$ and outward unit normal to the domain's boundary $\vec{n} = (n_x, n_y)$.

The goal here is to transform the above PDEs into a system of algebraic equations. To do so, equations~\eqref{eq:navier_2d} and~\eqref{eq:trac2} are considered in $N$ nodes, where corresponding partial differential operators are discretised using~\eqref{eq:ansatz}, resulting in a global system schematically represented as~\cite{SLAK20193}
\begin{equation}
	\label{eq:global-system}
	\begin{bmatrix}
		\mathbf{W}_{11} & \mathbf{W}_{12} \\  \mathbf{W}_{21} & \mathbf{W}_{22}
	\end{bmatrix}
	\begin{bmatrix}
		\mathbf{u} \\ \mathbf{v}
	\end{bmatrix} =
	\begin{bmatrix}
		\mathbf{b}_1 \\ \mathbf{b}_2
	\end{bmatrix},
\end{equation}
where $\mathbf{u}$ and $\mathbf{v}$ stand for vector of size $N$ containing the unknown discrete displacements components. Vectors $\b{b_1}$ and $\b{b_2}$ hold boundary conditions values and in blocks $\mathbf{W}_{11}, \mathbf{W}_{12}, \mathbf{W}_{21}, \mathbf{W}_{22}$ stencil weights~\eqref{eq:rbf-system-aug} are stored. To discretise Navier-Cauchy equation, stencil weights for approximation of first and second derivatives $\mathbf{w}_i^{\partial x}, \mathbf{w}_i^{\partial y}, \mathbf{w}_i^{\partial^2 x}, \mathbf{w}_i^{\partial^2 y}, \mathbf{w}_i^{\partial xy}$ are required in each node $i$ (how to compute those is discussed in section~\ref{sec:rbffd}). To assemble the system~\eqref{eq:global-system}, we need also a vector of stencil nodes of a $i$-th discretisation node $\mathcal{N}(i)$ that is in our case vector containing $n$ closest nodes. The system is finally expressed as

\begin{eqnarray}
	\left.
	\begin{aligned}
		W_{11}^{i, \mathcal{N}(i)_j} & = \left[ (\lambda + 2\mu)\mathbf{w}^{\partial^2}_i + \mu\mathbf{w}^{\partial^2 y}_i \right]_j \label{eq:u1} \\
		W_{12}^{i, \mathcal{N}(i)_j} & = \left[ (\lambda + \mu)\mathbf{w}^{\partial x\partial y}_i \right]_j                                       \\
		{\b{b}_1}_{i}                & = 0
	\end{aligned}
	\right\},
\end{eqnarray}
\begin{eqnarray}
	\left.
	\begin{aligned}
		W_{21}^{i, \mathcal{N}(i)_j} & = \left[ (\lambda + \mu)\mathbf{w}^{\partial x\partial y}_i \right]_j \label{eq:u2}           \\
		W_{22}^{i, \mathcal{N}(i)_j} & = \left[ \mu\mathbf{w}^{\partial^2 x}_i + (\lambda + 2\mu)\b{\chi}^{\partial^2 y}_i \right]_j \\
		{\b{b}_2}_{i}                & = 0                                                                                           \\
	\end{aligned}
	\right\},
\end{eqnarray}
with $i\in[1,N]$ and $i\in[1,n]$. In the same manner the lines representing the normal boundary conditions are assembled
\begin{eqnarray}
	\left.
	\begin{aligned}
		W_{11}^{i, \mathcal{N}(i)_j} & = \left[ \mu n_2 \mathbf{w}^{\partial y}_i + (2\mu + \lambda)n_1 \mathbf{w}^{\partial x}_i \right]_j \label{eq:t1} \\
		W_{12}^{i, \mathcal{N}(i)_j} & = \left[ \lambda n_1\mathbf{w}^{\partial y}_i +  \mu n_2 \mathbf{w}^{\partial x}_i \right]_j                       \\
		{\b{b}_1}_{i}                & = t_{0x,1}
	\end{aligned}
	\right\},
\end{eqnarray}
\begin{eqnarray}
	\left.
	\begin{aligned}
		W_{21}^{i, \mathcal{N}(i)_j} & = \left[ \mu n_1 \mathbf{w}^{\partial y}_i + \lambda n_2 \mathbf{w}^{\partial x}_i \right]_j \label{eq:t2} \\
		W_{22}^{i, \mathcal{N}(i)_j} & = \left[ \mu n_1\mathbf{w}^{\partial x}_i + (2\mu + \lambda)n_2 \mathbf{w}^{\partial y}_i \right]_j        \\
		{\b{b}_2}_{i}                & = t_{0y,i}
	\end{aligned}
	\right\},
\end{eqnarray}
while in essential boundary conditions nodes, we simply express
\begin{eqnarray}
	\begin{aligned}
		W11_{i, i}    & =1         \\
		{\b{b}_1}_{i} & = u_{0x,i}
	\end{aligned}
	& \quad \text{ and } \quad  &
	\begin{aligned}
		W22_{i, i}    & =1         \\
		{\b{b}_2}_{i} & = u_{0y,i}
	\end{aligned}
	\quad .
\end{eqnarray}
The example of assembled system is depicted in Figure~\ref{fig:spectra}. Solving this system provides the numerical solution to the Navier-Cauchy problem as formulated in Equation~\eqref{eq: navier-cauchy}.

\subsection{Note on implementation and used parameters}
The entire solution procedure is implemented in a C++~\footnote{Source code is available at \url{https://gitlab.com/e62Lab/public/p_2d_plasticity} under tag \emph{v1.0}.\label{fn:git}} and it is strongly dependent on our in-house a open-source mesh-free project \emph{Medusa}~\cite{slak2021medusa}. The code was compiled using \texttt{g++ (GCC) 11.3.0 for Linux} with \texttt{-O3 -DNDEBUG} flags on \texttt{AMD EPYC 7713P 64-Core Processor} computer. In cases where we show wall-clock times, the CPU frequency was fixed to 2.27 GHz, boost functionality disabled and CPU affinity assured using the \texttt{taskset} command. The highest RAM usage during the computation was approximately $\unit[4]{GB}$ out of available $\unit[12]{GB}$, thus, no RAM-related issues were encountered due to using the direct SparseLU solver in the process.

In all following analyses, a RBF-FD aproximation with PHS basis of order $k=3$ and augmenting monomials of order $m=3$ on stencil size of $n=20$ is used.
The Picard iteration tolerance is set to $\epsilon_{\mathrm{tol}} = 10^{-6}$
The final sparse linear system, is solved with the direct SparseLU using Eigen C++ library~\cite{eigenweb}.
Post-processing was done with Python 3.10.6 and Jupyter notebooks, also available in the provided git repository\footref{fn:git}.

\section{Numerical Examples}
\label{sec:examples}
The proposed solution procedure is applied to three different cases:
\begin{enumerate}[(i)]
	\item the linear-elastic case in Section~\ref{sec:first_example},
	\item the perfectly plastic elasto-plastic case in Section~\ref{sec:second_example} and
	\item the elasto-plastic case with linear isotropic hardening in Section~\ref{sec:isotropic}.
\end{enumerate}

The first two problems have been extensively studied by other researchers. De Souza et al.~\cite{desouza2008} used FEM to obtain numerical solutions while Hill~\cite{hill1998} derived the analytic solutions. Both analytical and FEM solutions are used here to evaluate the accuracy of present solution procedure. 

\begin{figure}[H]
	\centering
	\includegraphics[trim={17cm 17cm 0cm 0cm},clip, width=0.5\columnwidth]{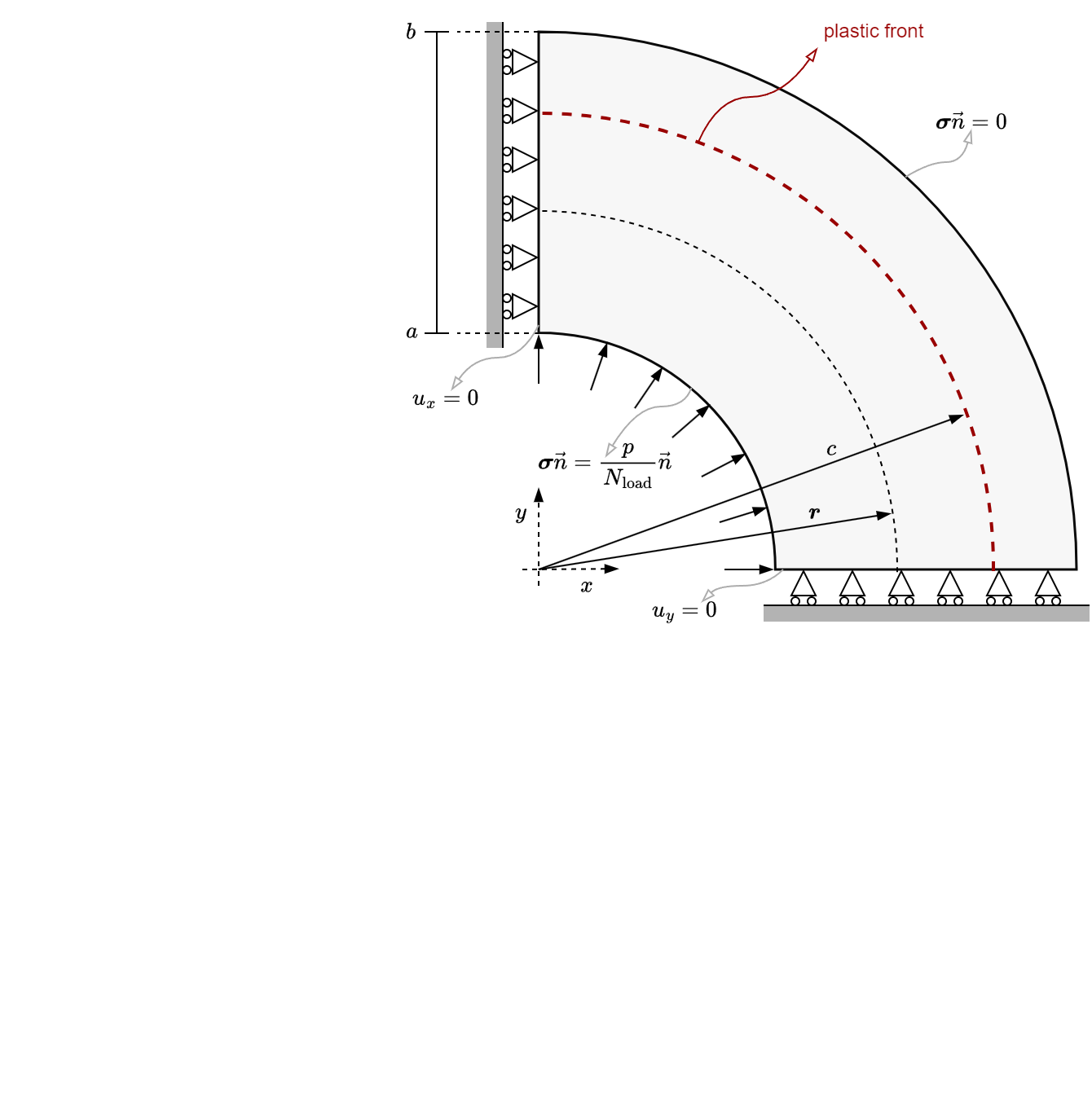}
	\caption{Schematic representation of a pressurized thick-walled pipe expansion.}
	\label{fig: sketch 2}
\end{figure}

The computational domain is the same for all three cases: a thick-walled pipe with inner and outer radii $a=\unit[100]{mm}$ and $b=\unit[200]{mm}$ respectively. As schematically presented in Figure~\ref{fig: sketch 2}, only a quarter of the pipe from the first quadrant is considered. The material properties were constant for all considered cases: Young's modulus $E=\unit[210]{GPa}$, Poisson's ratio $\nu = 0.3$ and yield stress $\sigma_y = \unit[0.24]{GPa}$.

The left side of the pipe is constrained to move freely along the $y$-direction, while the bottom side of the pipe is constrained to move only in the $x$-direction. The outer wall of the pipe is traction free, while the inner wall is subjected to internal pressure $p$, which, as explained in Section~\ref{sec:plasticity}, is gradually increased until the full magnitude is applied to the pipe's inner surface, i,e, at each load step $p/N_{\mathrm{load}}$ is applied. One can find all boundary conditions also on a problem scheme in Figure~\ref{fig: sketch 2}. 

The amount of internal pressure depends on the case being considered: $\unit[0.05]{GPa}$ for the linear-elastic case, $\unit[0.19]{GPa}$ for perfectly plastic case and $\unit[0.175]{GPa}$ for the case with linear isotropic hardening. Note that all of the following analyses are carried out assuming plane strain conditions, which means that the out-of-plane or $z$-direction strains are zero, i.e.\ $\ezz = \exz = \eyz = 0$.

\subsection{The Linear Elastic Case}
\label{sec:first_example}
\begin{figure}[h]
	\centering
	\includegraphics[width=\columnwidth]{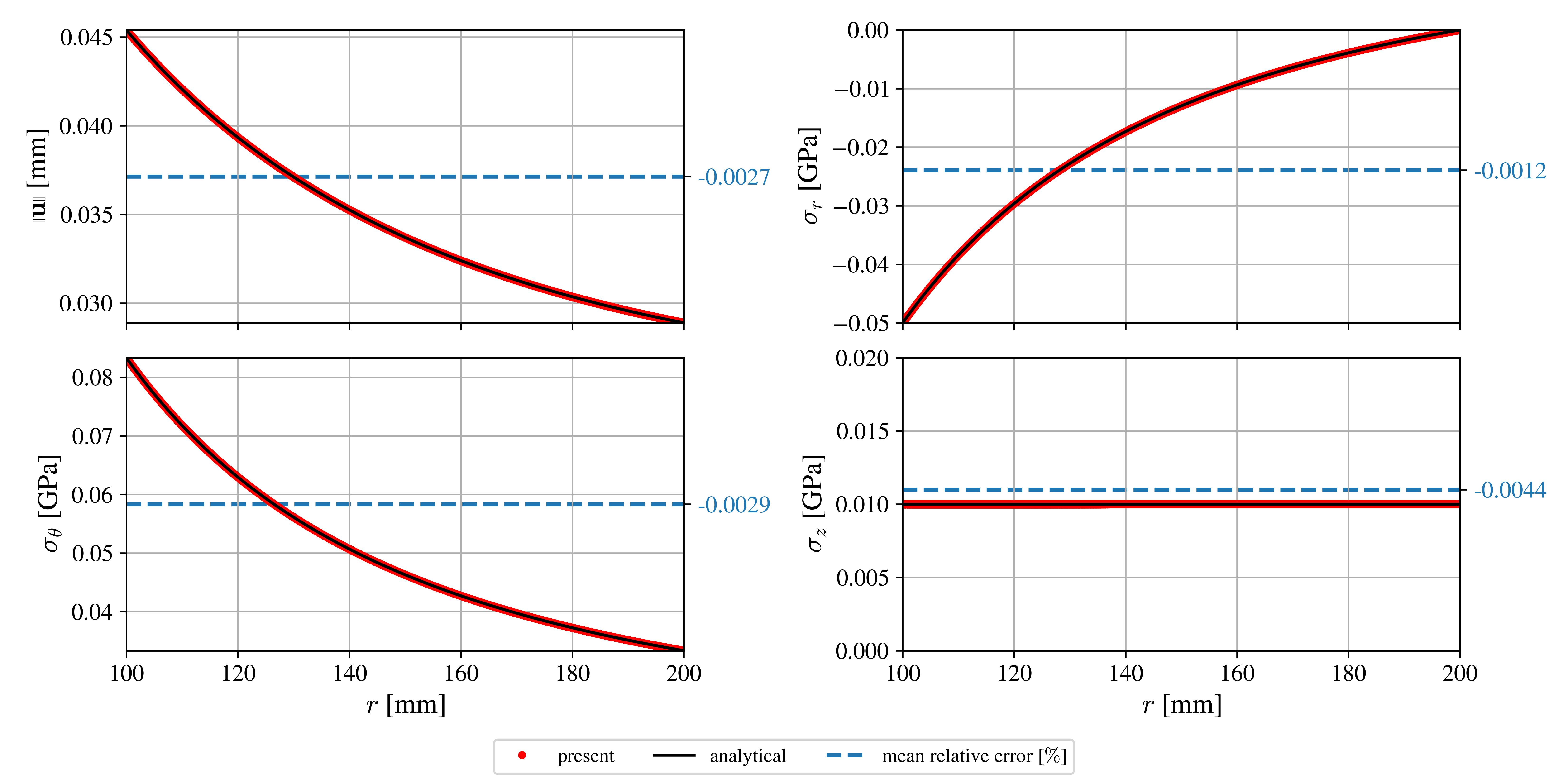}
	\caption{Assessment of the present solution of the purely elastic case. The scattered plot is displaying the relative error in all nodes. }
	\label{fig:elast_x_vs_r}
\end{figure}

\begin{figure}[H]
	\centering
	\includegraphics[width=0.285\textwidth]{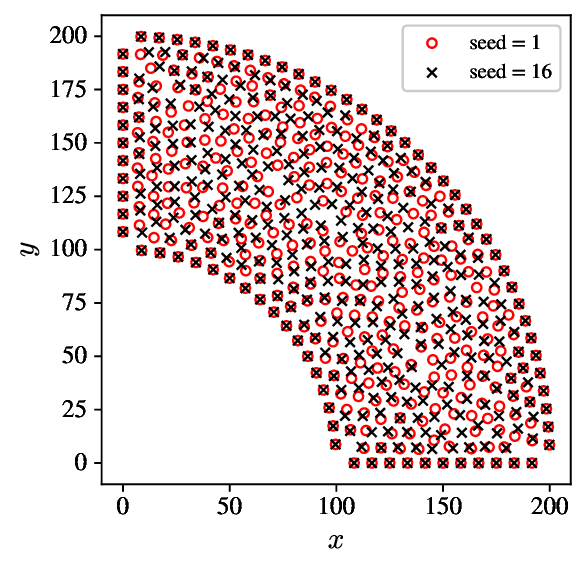}
	\includegraphics[width=0.33\linewidth]{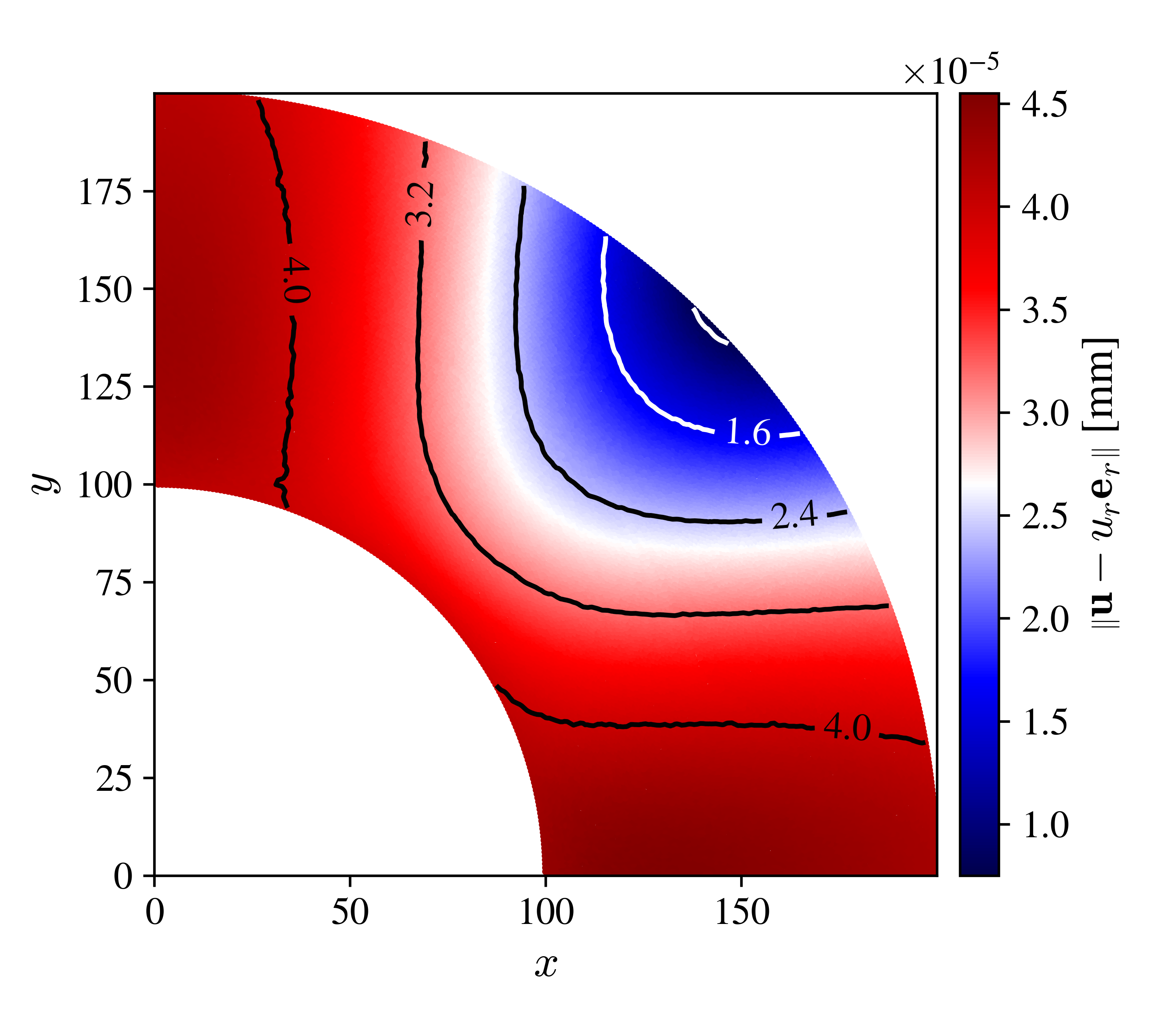}
	\includegraphics[width=0.325\linewidth]{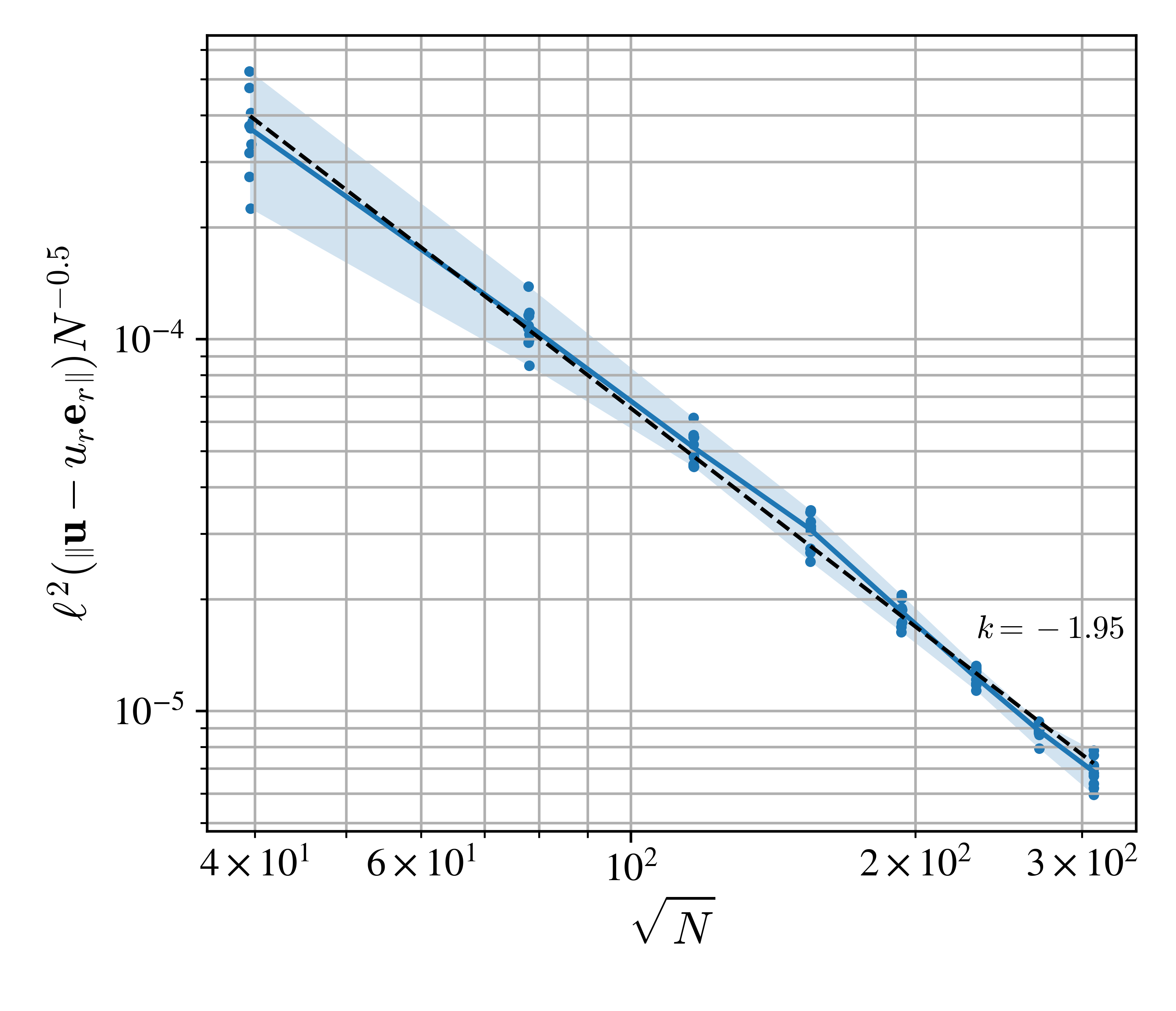}
	\caption{An example of  scattered nodes using two different random seeds (left), spatial distribution of the displacement magnitude error (middle) and convergence of the $\ell^2\mathrm{-norm}$ of the displacement vector for the linear-elastic case (right).}
	\label{fig:elast_disp_conv}
\end{figure}

To establish the basic confidence level in the proposed solution procedure, the internal pressure $p$ is set to a value of $\unit[0.05]{GPa}$, which is intentionally low, to assure that all deformations are well within the elastic regime. The case can be solved in a closed form~\cite{desouza2008,hill1998}
\begin{equation}
	\label{eq: analytical l-elastic}
	\begin{aligned}[c]
		\sigma_r        & = -\left(\frac{\frac{b^2}{r^2} - 1}{\frac{b^2}{a^2} - 1}\right) p, \\
		\sigma_{\theta} & = \left(\frac{\frac{b^2}{r^2} + 1}{\frac{b^2}{a^2} - 1}\right) p,
	\end{aligned}
	\qquad\qquad
	\begin{aligned}[c]
		\sigma_z & = \frac{2 \nu }{\frac{b^2}{a^2} - 1} p,                                                                                    \\
		u_r      & = \frac{p}{e} \frac{\left(1 + \nu\right)\left(1 - 2\nu\right) r + \left(1 + \nu\right) \frac{b^2}{r}}{\frac{b^2}{a^2} - 1}
	\end{aligned}
\end{equation}
in cylindrical coordinates, which we will use for a proper verification of the Navier-Cauchy solver. To verify the performance of the proposed solution procedure using scattered nodes, we numerically solve the problem in Cartesian coordinate system (left plot of Figure~\ref{fig:elast_disp_conv}). In Figure~\ref{fig:elast_x_vs_r}) the RBF-FD solution is compared with the~\eqref{eq: analytical l-elastic}.
We can see that the present method is in good agreement with the closed-form solution, with the relative error of the numerical solution well below 1 $\permil$ using $N\approx95\,000$ nodes.
In right plot of Figure~\ref{fig:elast_disp_conv}, we vary the internodal distance $h$, examining its effect on the numerical solution's accuracy as the local field description improves. The error $\norm{{\mathbf u}( \x_i) - u_r(\x_i)\mathbf e_r(\x_i)}$ is computed for each computational point $\mathbf x_i \in \Omega$ and the $\ell^2$-norm of the numerical solution error is plotted against the number of discretization nodes $N$. For clarity, an example spatial distribution of the displacement magnitude error is shown in the middle of Figure~\ref{fig:elast_disp_conv}.
Since there is practically an infinite number of possible node layouts for a given $h$, depending on the number and position of the seed nodes, the number of expansion candidates and different randomizations (random \emph{seed}) in candidate generation, we also test the spread of the solutions in relation to this variation. The left plot of Figure~\ref{fig:elast_disp_conv} shows two different node layouts with different initial node generation settings. For each data point in the right plot of Figure~\ref{fig:elast_disp_conv}, we varied the node generation setting and obtained an ensemble of $10$ solutions as a result. The ensemble converges in both respects, i.e. the finer the domain discretizations the higher the accuracy and the lower the spread.

Note that all solutions in this section have been obtained by immediately enforcing the full amount of internal pressure $p$, omitting the gradual load increase over $\Nload$ load steps for computation by setting $\Nload = 1$. This is justified as only elastic deformations are expected and the solution procedure never encountered a local violation of the von Mises yield criterion from Equation~\eqref{eq: y funct}, i.e.\ the return mapping Algorithm~\ref{alg:return_mapping} was never used.

\begin{figure}[h]
	\centering
	\includegraphics[width=\linewidth]{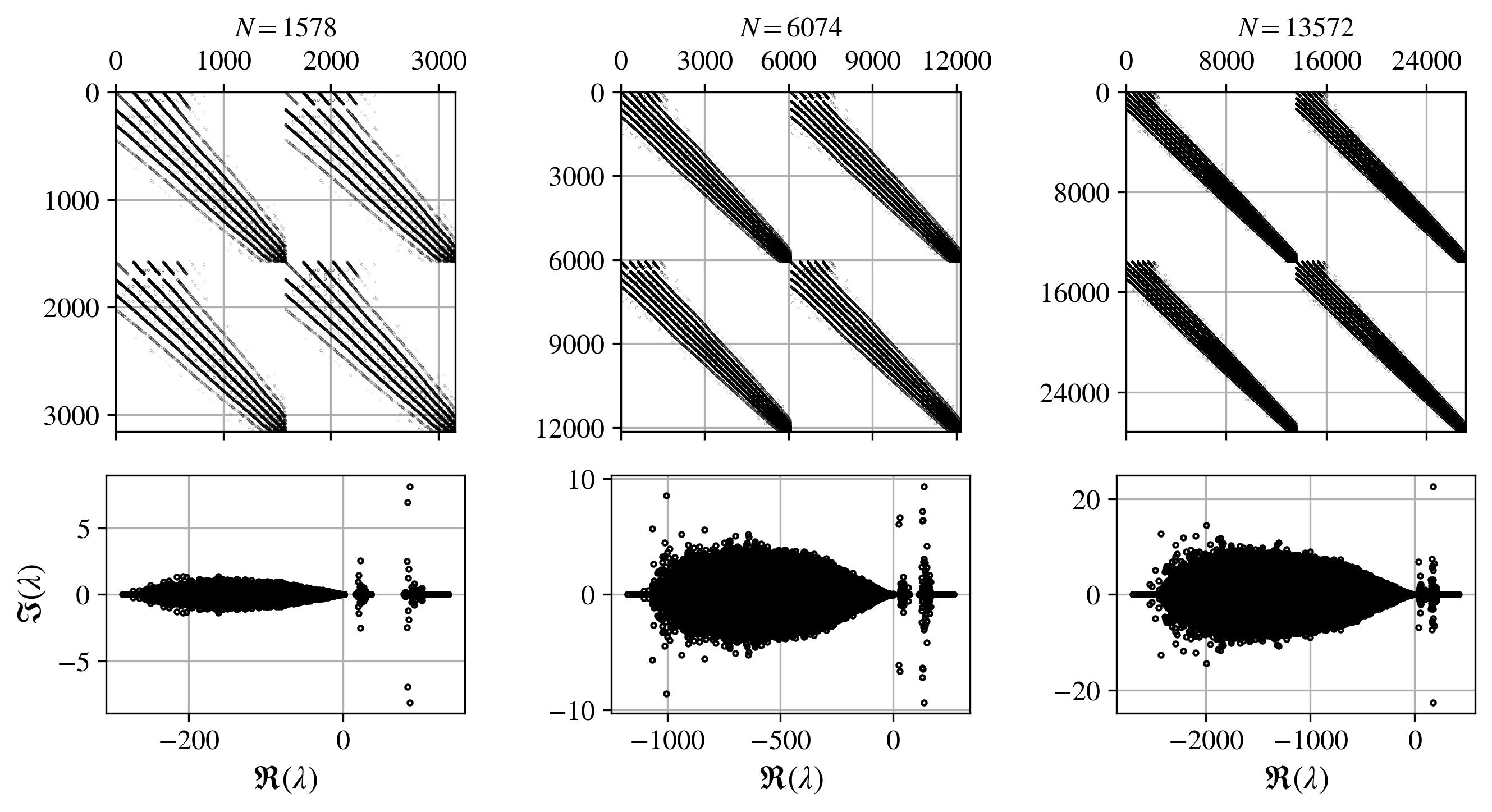}
	\caption{Sparsity of final matrices (top row) and spectrum (bottom row) for three different discretization qualities.}
	\label{fig:spectra}
\end{figure}

The sparsity patterns of the final large sparse systems are shown in Figure~\ref{fig:spectra}. The top row shows the sparsity of the systems for three different domain discretization qualities. The bottom row of Figure~\ref{fig:spectra} shows the spectra of the three matrices. The ratios between the real and imaginary parts of the eigenvalues are in good agreement with previous studies~\cite{janvcivc2021monomial,bayona2017role,bayona2019comparison}. Most of the eigenvalues have an order of magnitude larger real-valued part compared to the relatively small imaginary parts.

\subsection{Perfectly-Plastic Yielding of an Internally Pressurized Thick-Walled Cylinder}
\label{sec:second_example}

With the established confidence in the implementation of elastic deformation, the internal pressure $p$ is increased to $\unit[0.19]{GPa}$, which is enough to introduce the irreversible plastic deformations. The pressure $p$ is gradually increased until a prescribed load is reached, as proposed in Algorithm~\ref{alg:adapt}.

For the present cylindrically symmetrical problem, the plastic yielding starts at the inner surface ($r=a$) and gradually develops towards the outer surface ($r=b$) in the form of a circular front located at $r=c$. The material collapses when the plastic front reaches the outer surface ($c=b$). At that point, the entire domain is plastically deformed and the external load required to achieve that state is referred to as the \emph{limit load}\footnote{For the present dimensions and material properties, the limit load is $p_\mathrm{lim}\approx \unit[0.19209]{GPa}$ as discussed by de Souza~\cite{desouza2008}.}. 

Analytical solutions for the normal stress tensor components and for the radial displacement are given in~\cite{hill1998,chakrabarty1987}. Hill, in his derivation, assumes the Tresca yielding criterion, but argues that this can make a fair approximation for the von Mises yielding criterion~\cite{hill1998}. Provided that the cylinder has started yielding, and the position of the plastic front within the cylinder is located at $c$, the analytical approximations read

\begin{equation}
	\label{eq: analytical tresca}
	\begin{aligned}
		    & \begin{aligned}[c]
			      \left.
			      \begin{aligned}
				\sigma_r        & = -\frac{2  c^2}{\sqrt{3} b^2} \left(\frac{b^2}{r^2} - 1\right)\sigma_y \\
				\sigma_{\theta} & = \frac{2 c^2}{\sqrt{3} b^2} \left(\frac{b^2}{r^2} + 1\right) \sigma_y  \\
				\sigma_z        & = \frac{4 \nu  c^2}{\sqrt{3} b^2}\sigma_y
			\end{aligned}
			      \right\}
			      c \leq r \leq b, \\
		      \end{aligned}
		\qquad\qquad
		\begin{aligned}[c]
			\left.
			\begin{aligned}
				\sigma_r        & = -\frac{2 }{\sqrt{3}} \left(1 - \frac{c^2}{b^2} + \ln \frac{c^2}{r^2}\right) \sigma_y   \\
				\sigma_{\theta} & = \frac{2}{\sqrt{3}} \left(1 + \frac{c^2}{b^2} + \ln \frac{c^2}{r^2}\right) \sigma_y     \\
				\sigma_z        & = \frac{4 \nu  c^2}{\sqrt{3}} \left(\frac{c^2}{b^2} - \ln \frac{c^2}{r^2}\right)\sigma_y
			\end{aligned}
			\right\}
			a \leq r \leq c,
		\end{aligned}
		\\
		u_r & = \left(1 - \nu\right) \frac{2 \sigma_y c^2}{\sqrt{3} \mu r} + \left(1 - 2 \nu\right) \frac{\sigma_r r}{2 \mu}, \quad a \leq r \leq b,
	\end{aligned}
\end{equation}
where $c$ is determined from its relation to the cylinder's internal pressure $p$ as
\begin{equation}
	\label{eq: c determine}
	p = \frac{2 \sigma_y}{\sqrt{3}} \left(1 - \frac{c^2}{b^2} + \ln \frac{c^2}{a^2}\right).
\end{equation}

\begin{figure}[H]
	\centering
	\includegraphics[width=0.86\columnwidth]{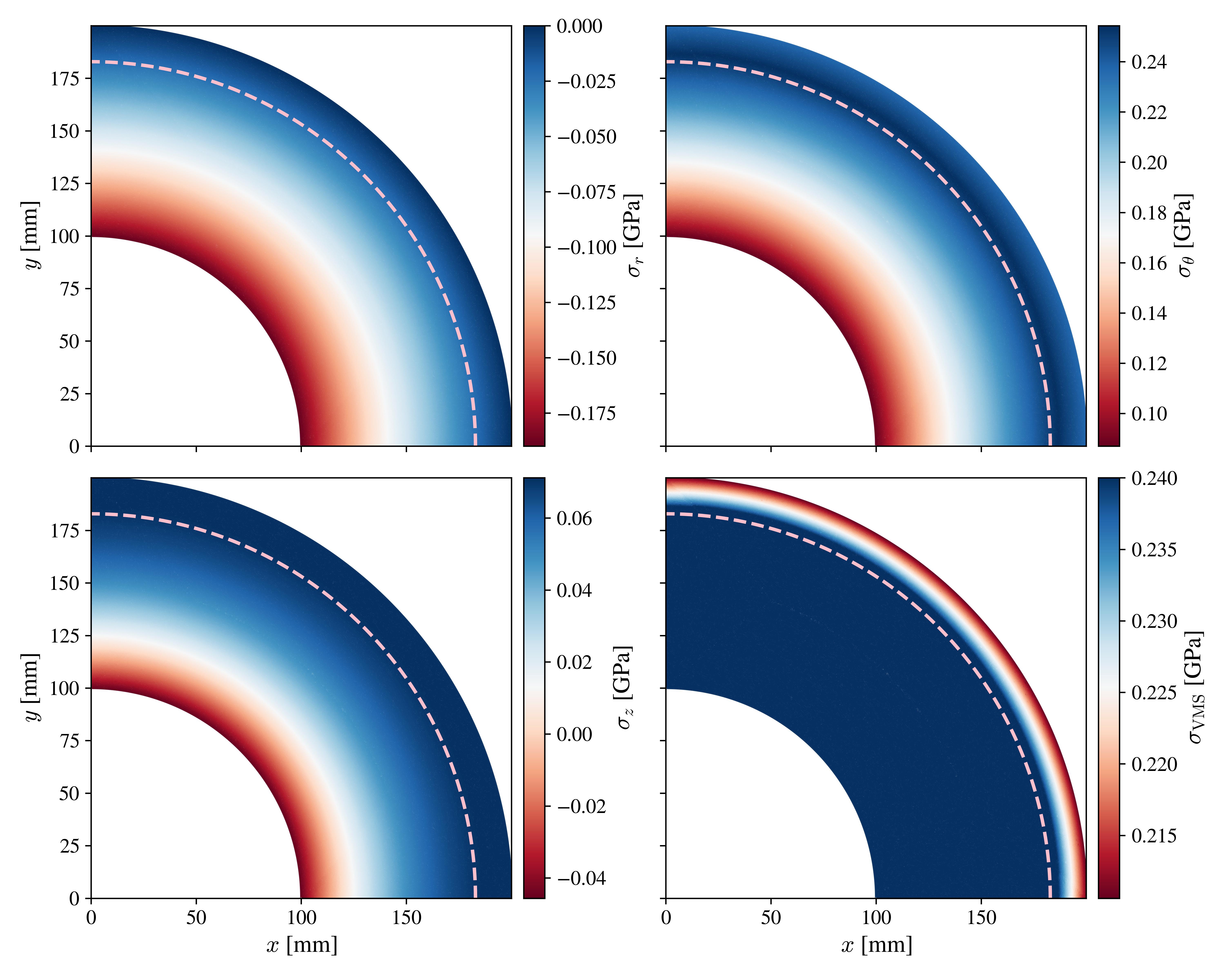}
	\caption{Example of numerical solution for perfectly plastic case at $p = \unit[0.19]{GPa}$, $N=95\,300$ nodes at the final load step $\Nload=25$. The dashed line represents plastic front position $c=186.67 \ \mathrm{mm}$.}
	\label{fig: eg result}
\end{figure}

As in the previous case, we use Cartesian coordinates, even though the solution is given in cylindrical coordinates. In this case, the internal pressure has been gradually increased over $\Nload = 25$ load steps. Example solution for the last load step corresponding to internal pressure $p= \unit[0.19]{GPa}$ is shown in Figure~\ref{fig: eg result}. Figure shows the normal stresses $\sigma_r$, $\sigma_\theta$ and $\sigma_z$ and the von Mises stress $\svms$. Additionally, the plastic front at $c=186.67 \ \mathrm{mm}$ is shown with a dashed pink line.

\begin{wrapfigure}{r}{0.5\columnwidth}
	\vspace{-0.5cm}
	\centering
	\includegraphics[width=0.5\columnwidth]{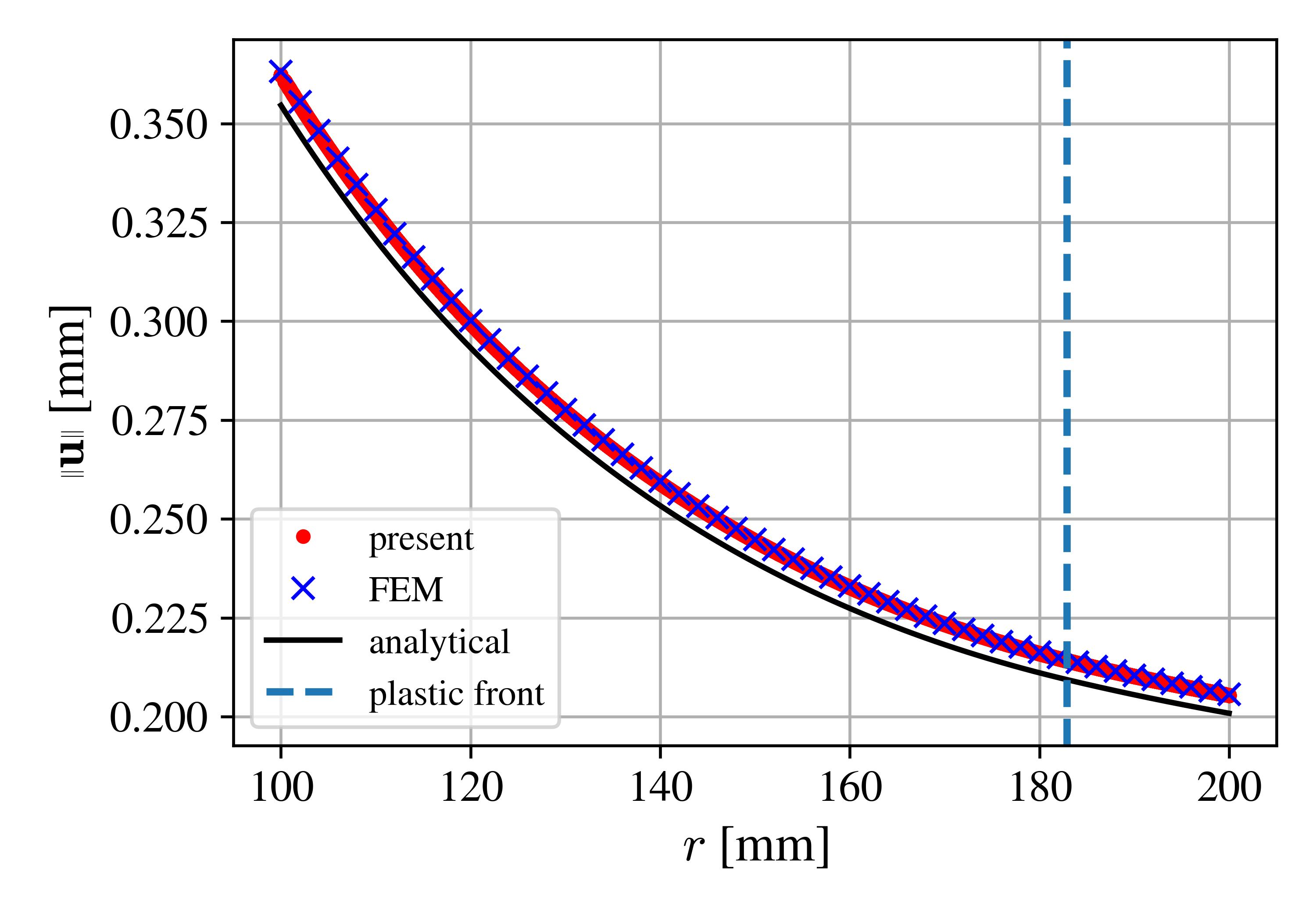}
	\caption{Displacement magnitude $\norm{\mathbf u}$ for the perfectly plastic case. For clarity, only every second point is shown for the FEM solution.}
	\label{fig: u vs fem}
\end{wrapfigure}
\begin{figure}[H]
	\centering
	\includegraphics[width=0.9\columnwidth]{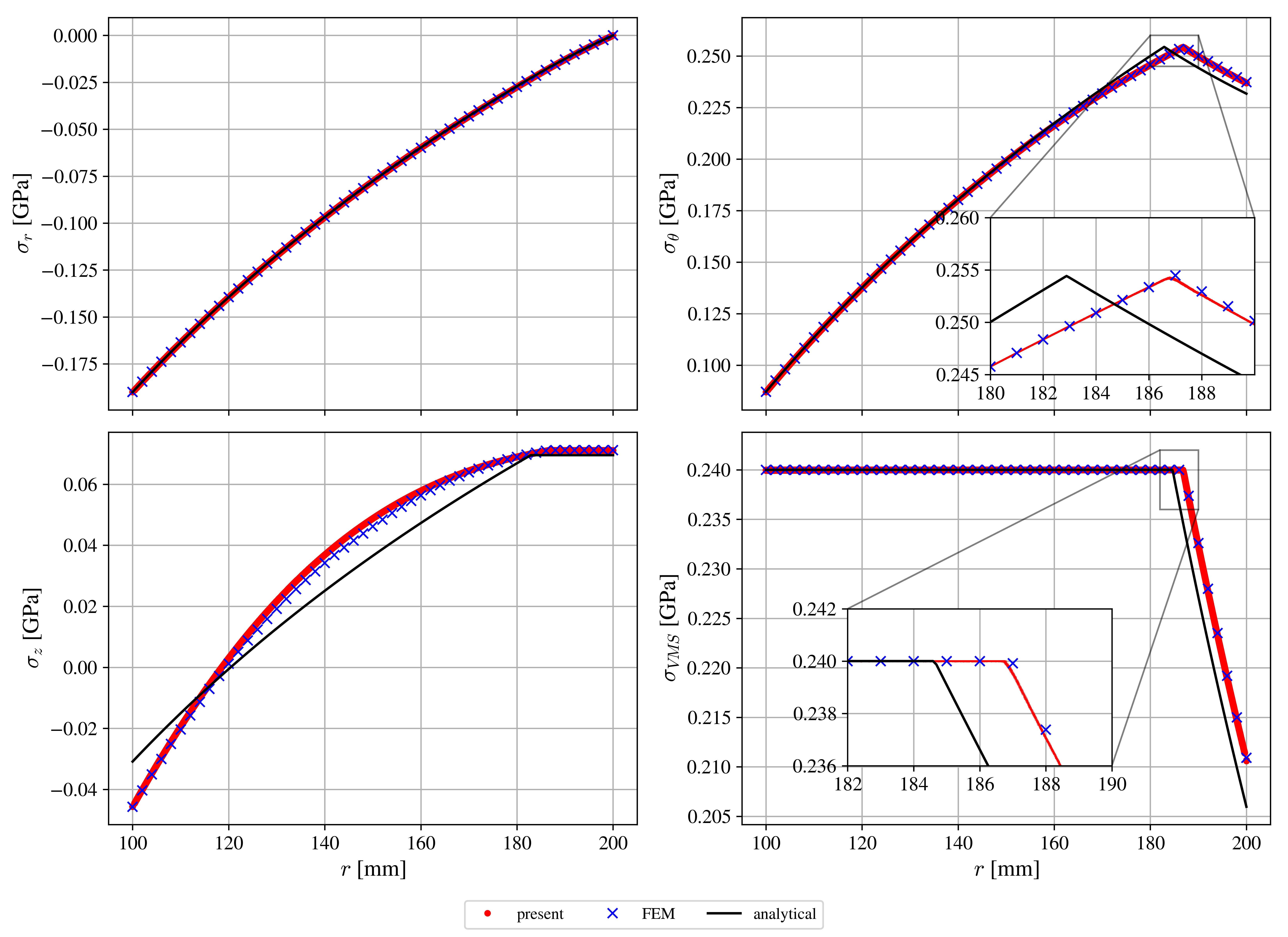}
	\caption{Stress analysis in case of perfectly elasto-plastic problem.}
	\label{fig:perfect_stress}
\end{figure}

The accuracy of the numerical solution is further evaluated in Figures~\ref{fig: u vs fem} and~\ref{fig:perfect_stress}, where the present solution procedure is compared to the analytic solution (see Equations~\eqref{eq: analytical tresca}) and to the FEM solution. To obtain the FEM solution, commercial software package Abaqus was used. The FEM domain was discretized with a structured mesh consisting of 23\,600 8-node biquadratic plane strain quadrilaterals (\texttt{CPE8}) with $71\,473$ nodes, corresponding to element side length of approximately $\unit[1]{mm}$.

In Figure~\ref{fig: u vs fem} the displacement magnitude with respect to the radial coordinate $r$ is shown for all $N=95\,300$ nodes of the present method and for nodes along the $y=0$ line for the FEM approximation. We are pleased to observe that the proposed solution procedure is in good agreement with the FEM solution throughout the entire domain. Both FEM and RBF-FD solutions, however, are slightly offset from the analytical solution. The most likely reason for the mismatch compared to the analytical solution is the different yield criteria used.

Moreover, in Figure~\ref{fig:perfect_stress}, the performance of the present method is evaluated in terms of stresses. Specifically, the figure shows the normal stress components and the von Mises stress $\svms$ with respect to the radial coordinate $r$. Good agreement of the present method and the FEM solution is observed, while both again show some differences compared to the analytic solution, particularly near the elastic-plastic interface.

\subsubsection{The effect of load stepping}
\label{sec:load_steps}

\begin{figure}[H]
	\centering
	\includegraphics[width=0.45\linewidth]{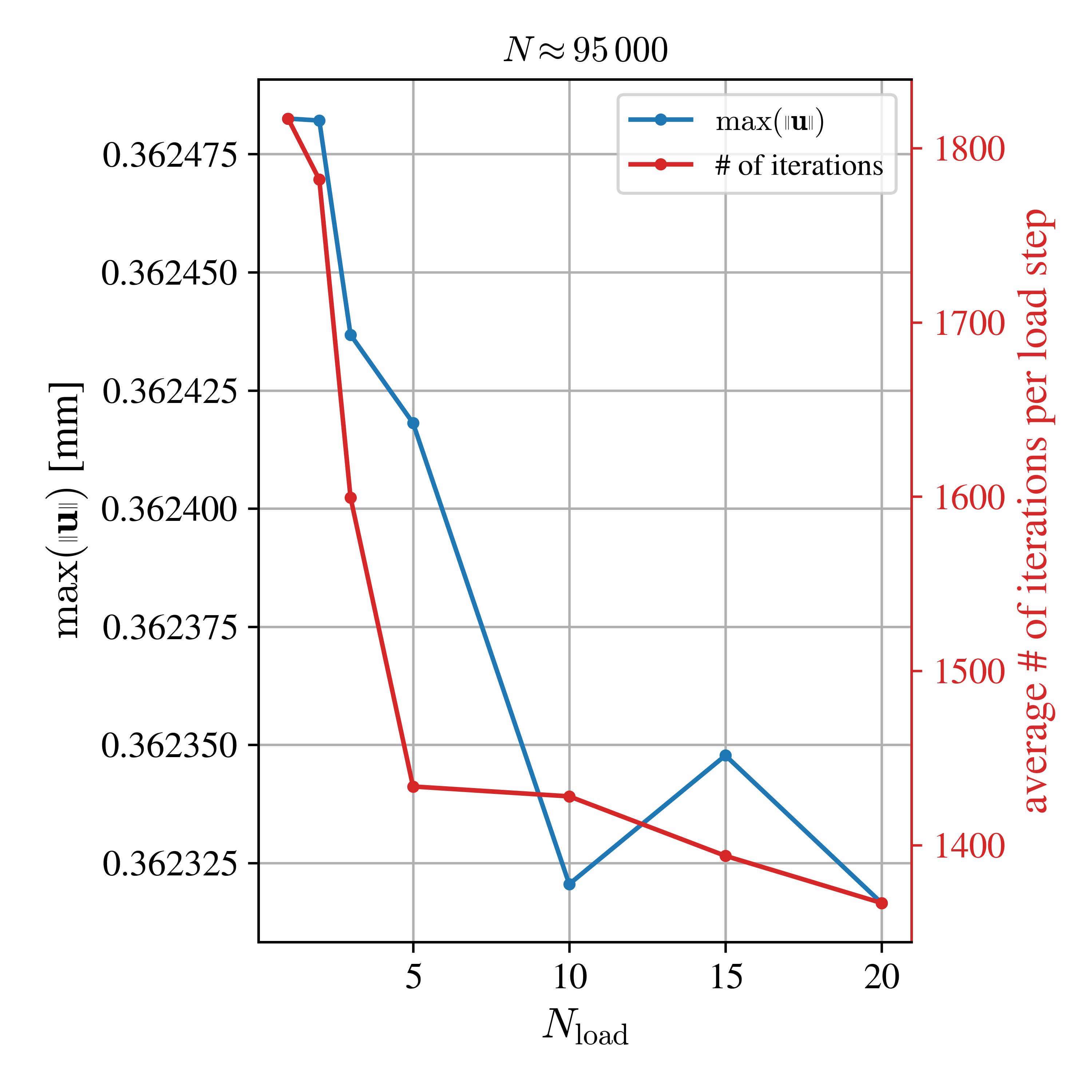}
	\includegraphics[width=0.45\linewidth]{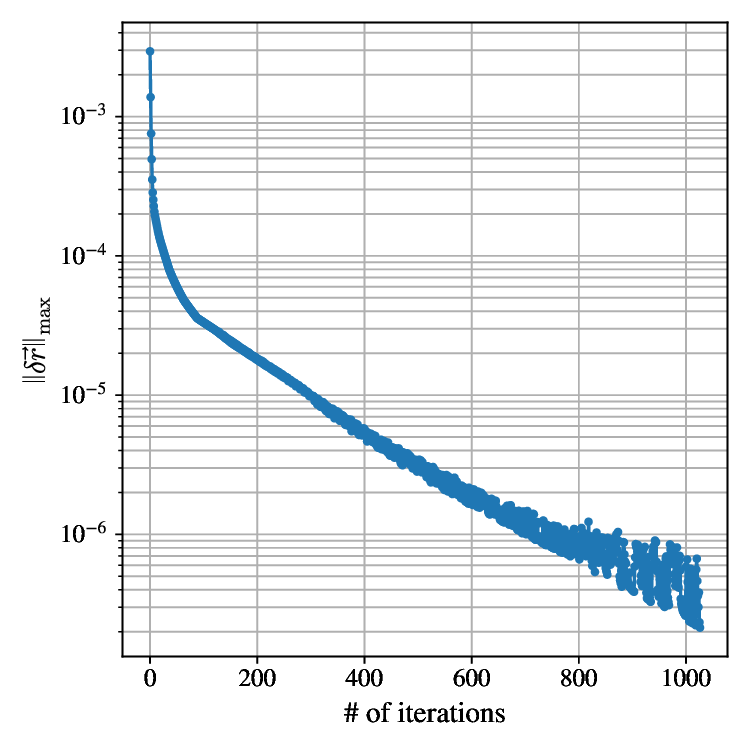}	
	\caption{The effect of number of load steps $\Nload$ on the accuracy of the numerical solution in terms of maximum displacement magnitude $\norm{\mathbf u}$ and number of Picard iteration per step (left). Example of residual convergence at sixth load step in simulation with ten load steps (right).}
	\label{fig:load_analysis}
\end{figure}

This section briefly discusses the impact of  $\Nload$ on solution accuracy and the  convergence of the Picard iteration. In the left plot of Figure~\ref{fig:load_analysis}, an effect of $\Nload$ on the final displacement accompanied with the average number of Picard iterations per step is presented. We see that in this particular case $\Nload$ has no effect on the final solution, i.e. the maximum displacement is almost constant with respect to the $\Nload$ (note the scale of the left $y$-axis). On the other hand, increasing the $\Nload$ leads to additional computational effort. Although the average number of required Picard iterations per step decreases with increasing $\Nload$, the cumulative number of iterations increases. The right plot of Figure~\ref{fig:load_analysis} shows the Picard convergence for a setup with $\Nload=10$ in the sixth step, which is similar to the one presented in~\cite{yarushina2010}.  
 
However, it is important to consider that, for some problems, the degree of non-linearity of the return-mapping equations is so high that the convergence radius of iteration scheme becomes reduced~\cite{desouza2008}. In such cases, improved initial guesses or a increased number of load steps is recommended~\cite{desouza2008,krabbenhoft2002basic}. Some researchers also proposed more sophisticated algorithms~\cite{desouza2008,krabbenhoft2002basic,PEREZFOGUET20014627} in order to improve the convergence.

\subsubsection{The Plastic Front Position}
\label{sec:perfect_front}
The following set of analyses is devoted to study the behaviour of hoop ($\sigma_\theta$) and axial ($\sigma_z$) stresses in Figure~\ref{fig:perfect_stress}. Neither of the two numerical solutions is able to accurately follow the analytic solution. Nevertheless, a zoomed in section of the hoop stress shows the breaking point, positioned at the elastic-plastic interface, hints that the two conceptually different numerical solutions are in good agreement. It appears as if both numerical approaches overestimate the position of the plastic front by approximately $\unit[4]{mm}$ compared to the analytic solution (see also Table~\ref{tb: res perfect}). A possible reason for this observed discrepancy between analytic and numerical solutions may be in the derivation of the analytic solution, originally derived for the Tresca yielding criterion. The Tresca yielding theory\footnote{From the design point of view, the Tresca yield condition is safer, but it also leads the engineer to take unnecessary measures in attempt to prevent an unlikely failure. Experiments suggest that the von Mises yield criterion is the one which provides better agreement with observed behaviour than the Tresca yield criterion. However, the Tresca yield criterion is still used because of its mathematical simplicity.} is known to be more conservative than the von Mises theory~\cite{ROOSTAEI202223}.

\begin{figure}[H]
	\centering
	\includegraphics[width=0.9\columnwidth]{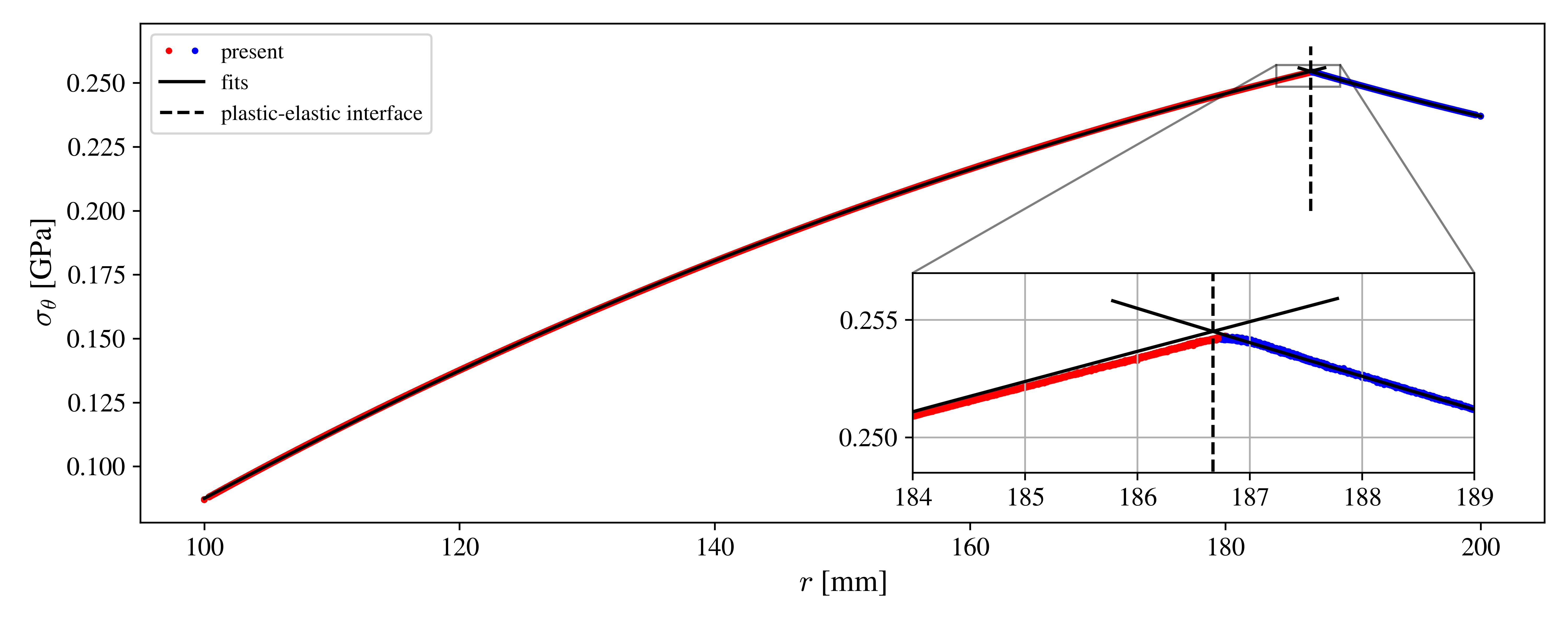}
	\caption{Example of plastic-elastic interface determination. Blue dots represent the points where the accumulated plastic strain $\bar{\varepsilon}^p = 0$, and red dots $\bar{\varepsilon}^p > 0$.}
	\label{fig: c determine}
\end{figure}

The position of plastic front $c$ is estimated from the numerical solution $\widehat{u}$. Example demonstrating how the elastic-plastic interface has been determined is shown in Figure~\ref{fig: c determine}. Red and blue points are used to mark the plastic and elastic regime, respectively. To capture the front position we fitted following function to the computed data
\begin{equation}
	f_\mathrm{fit}(r) = C_1 r^{-1} + C_2 r^{-2} + C_3 \log(r^{-1}) + C_4
\end{equation}
with constants $C_i$ (provided in Table~\ref{tab:fit}) to fit the data and a full black line to denote them in Figure~\ref{fig: c determine}.
The intersection of the two fitted lines, gives us an estimate on the position of the plastic front $c$. For clarity, only analysis for the present solution procedure is shown in Figure~\ref{fig: c determine} while Table~\ref{tb: res perfect} also shows the numerical data for the FEM solution.

\begin{table}[H]
	\centering
	\caption{Constant values to reproduce the fit lines from Figure~\ref{fig: c determine}.}
	\label{tab:fit}
	\def\arraystretch{1.3}
	\begin{tabular}{l|cccc}
		\toprule\toprule
		regime  & $C_1$                & $C_2$              & $C_3$                 & $C_4$    \\\midrule
		plastic & $-9.015\cdot 10^{1}$ & $2.751\cdot 10^3$  & $ 8.877\cdot 10^{-2}$ & $1.123$  \\
		elastic & $ 3.036\cdot 10^{2}$ & $-9.989\cdot 10^3$ & $-7.823\cdot 10^{-1}$ & $-5.176$ \\
		\bottomrule
	\end{tabular}
\end{table}

Table~\ref{tb: res perfect} clearly shows that the two numerical solutions are in good agreement in terms of all studied variables. The von Mises stress in both extremal parts of the domain, i.e.\ $r=a$ and $r=b$ are in good agreement. The displacement magnitude $\norm{\mathbf u}$ shows a marginal difference but the plastic front position $c$ again differs for only one tenth of a millimetre. As previously observed in Figure~\ref{fig:perfect_stress}, the closed-form solution predicts the elastic-plastic interface closer to the origin by approximately $\unit[4]{mm}$.

\begin{table}[H]
	\centering
	\caption{Results of the perfect plastic case - von Mises stress, and displacement magnitude at inner and outer surfaces, and position of the elastic-plastic interface. Under ``Present'' results we present the solution with $N \approx 95\,000$ vs.\ FEM vs.\ analytical approximation.}
	\label{tb: res perfect}
	\begin{threeparttable}[b]
		\def\arraystretch{1.3}
		\begin{tabular}{l|cc|cc|c|c}
			\toprule\toprule
			           & \multicolumn{2}{c}{$\svms$ [GPa]} & \multicolumn{2}{|c|}{$\norm{\mathbf u}$ [mm]} & \multirow{2}{*}{$c$ [mm]}  & \multirow{2}{*}{$N$}                                           \\
			           & $r = a$                           & $r = b$                                       & $r = a$                    & $r = b$                    &                                   \\
			\hline
			Present    & 0.2400\tnote{\mjstarsmall}        & 0.2106\tnote{\mjstarsmall}                    & 0.3623\tnote{\mjstarsmall} & 0.2053\tnote{\mjstarsmall} & 186.67\tnote{\mjstar} & $95\,300$ \\
			FEM        & 0.2400                            & 0.2109                                        & 0.3593                     & 0.2057                     & 186.79\tnote{\mjstar} & $71\,473$ \\
			Analytical & 0.2400\tnote{\mjtri}              & 0.2060                                        & 0.3546                     & 0.2008                     & 182.89                & -         \\
			\bottomrule
		\end{tabular}
		\begin{tablenotes}
			\item [\mjstarsmall] Average for all nodes on the edge.
			\item [\mjstar] Determined by fitting lines through points that are in the plastic zone ($\varepsilon^p > 0$), and points that are outside it ($\varepsilon^p = 0$), and finding the intersection between the two fits.
			\item [\mjtri] Presumed von Mises stress criterion.
		\end{tablenotes}
	\end{threeparttable}
\end{table}

\subsubsection{Shape of the Elastic-Plastic Front}
\label{sec:front_shape}
\begin{figure}[H]
	\centering
	\includegraphics[width=\columnwidth]{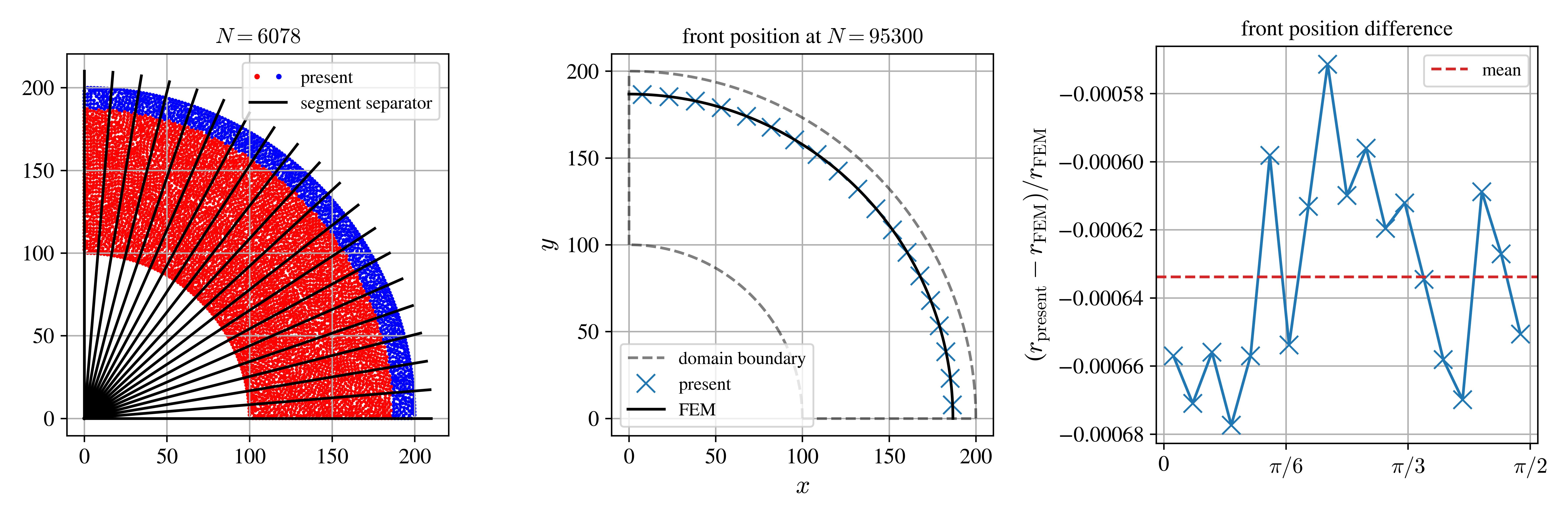}
	\caption{The elastic-plastic front distance to the origin computed for each of the $\Nseg=20$ segments.} 
	\label{fig:scheme_segments}
\end{figure}

Due to the cylindrical symmetry of the governing problem, the obtained numerical solution is expected to show symmetrical properties. For instance, the elastic-plastic front shape is expected to preserve a circular shape. Therefore, in this section, the front positioning algorithm from Section~\ref{sec:perfect_front} is extended and applied to $\Nseg$ segments, which equally divide the computational domain. A schematic presentation of the domain division is shown in Figure~\ref{fig:scheme_segments} on the right, where a coarse domain discretization with $N=6078$ is shown for clarity.

To determine the front position within each segment, the front position is estimated using the same procedure as in Section~\ref{sec:perfect_front}, except that only nodes within a given segment are used compared to the algorithm in Section~\ref{sec:perfect_front} working on all computational nodes.

The shape of the front has been computed for the finest domain discretization. Results for $\Nseg=20$ segments with approximately $5000$ nodes within each segment are shown in Figure~\ref{fig:scheme_segments} in the centre. For reference, the fully-symmetrical FEM solution is added with a full black line. The numerical solution obtained with the present approach indeed preserves the circular shape of the elastic-plastic front position. The circular shape is further evaluated in Figure~\ref{fig:scheme_segments} on the right, where the distance of the plastic front from the origin $r_\mathrm{present}$ is compared to distance $r_\mathrm{FEM}$ of the FEM-based solution.   

\subsection{Linearly-Hardening Yielding of an Internally Pressurized Thick-Walled Cylinder}
\label{sec:isotropic}
In this section, we consider the case of isotropic linear hardening in the plastic regime (see dotted red line in Figure~\ref{fig:stress_strain_sketch}). Although pure linear hardening is hardly observed in real materials, it forms strong foundations for the simulation of realistic materials. The hardening curve for real materials is typically measured experimentally during an uniaxial tension test, where the experimental data is interpolated with a piecewise-linear function to create a model that can be used in numerical simulation.
\begin{wrapfigure}{r}{0.5\columnwidth}
	\centering
	\includegraphics[width=0.5\columnwidth]{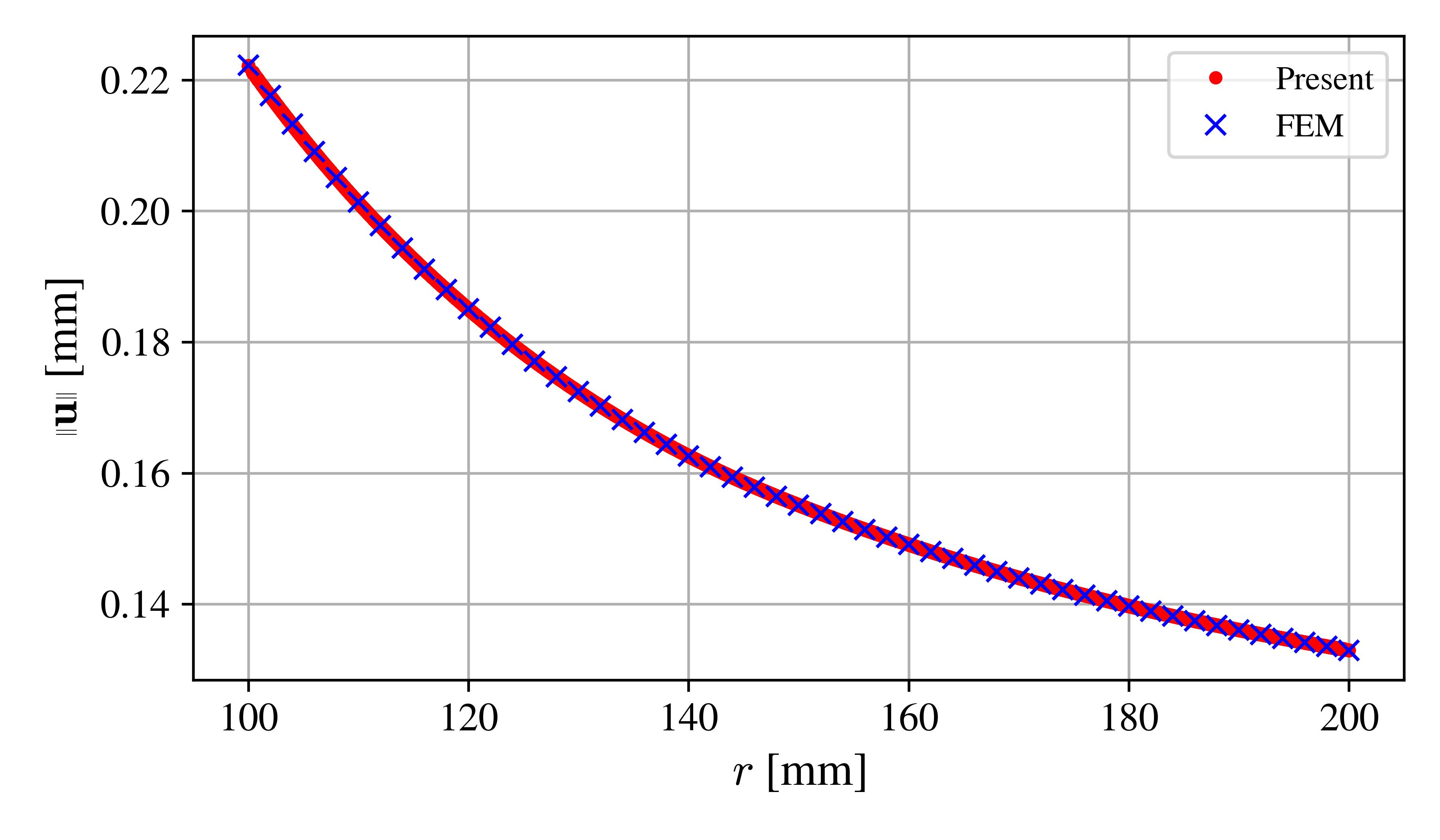}
	\caption{Displacement magnitude $\norm{\mathbf u}$ vs radial coordinate $r$. Present solution compared with FEM. Note that the FEM solution is only displaying every second point for clarity.}
	\label{fig: u vs fem lin}
\end{wrapfigure}
Specifically, in our case, the hardening curve is defined as follows
\begin{equation}
	\sigma_y(\overline \varepsilon^p) = \unit[0.24]{MPa} + \unit[10]{MPa} \cdot  \overline \varepsilon^p,
\end{equation}
resulting in a constant hardening slope $H = \unit[10]{MPa}$ throughout the entire plastic regime.

Example solution fields obtained with the present solution procedure are similar to the ones shown in Figure~\ref{fig: eg result} and are, thus, not shown again. Instead, the numerical solution obtained with the present method is compared to the community standard -- the FEM solution from Abaqus. The Abaqus setup was the same as in the previous subsection. The two solutions are compared in terms of the displacement field magnitude with respect to the radial coordinate in Figure~\ref{fig: u vs fem lin}. Good agreement between the two conceptually different approaches is observed throughout the entire domain in radial coordinates. These observations are further confirmed with numerical data in Table~\ref{tb: res linear}.

\begin{figure}[H]
	\centering
	\includegraphics[width=\columnwidth]{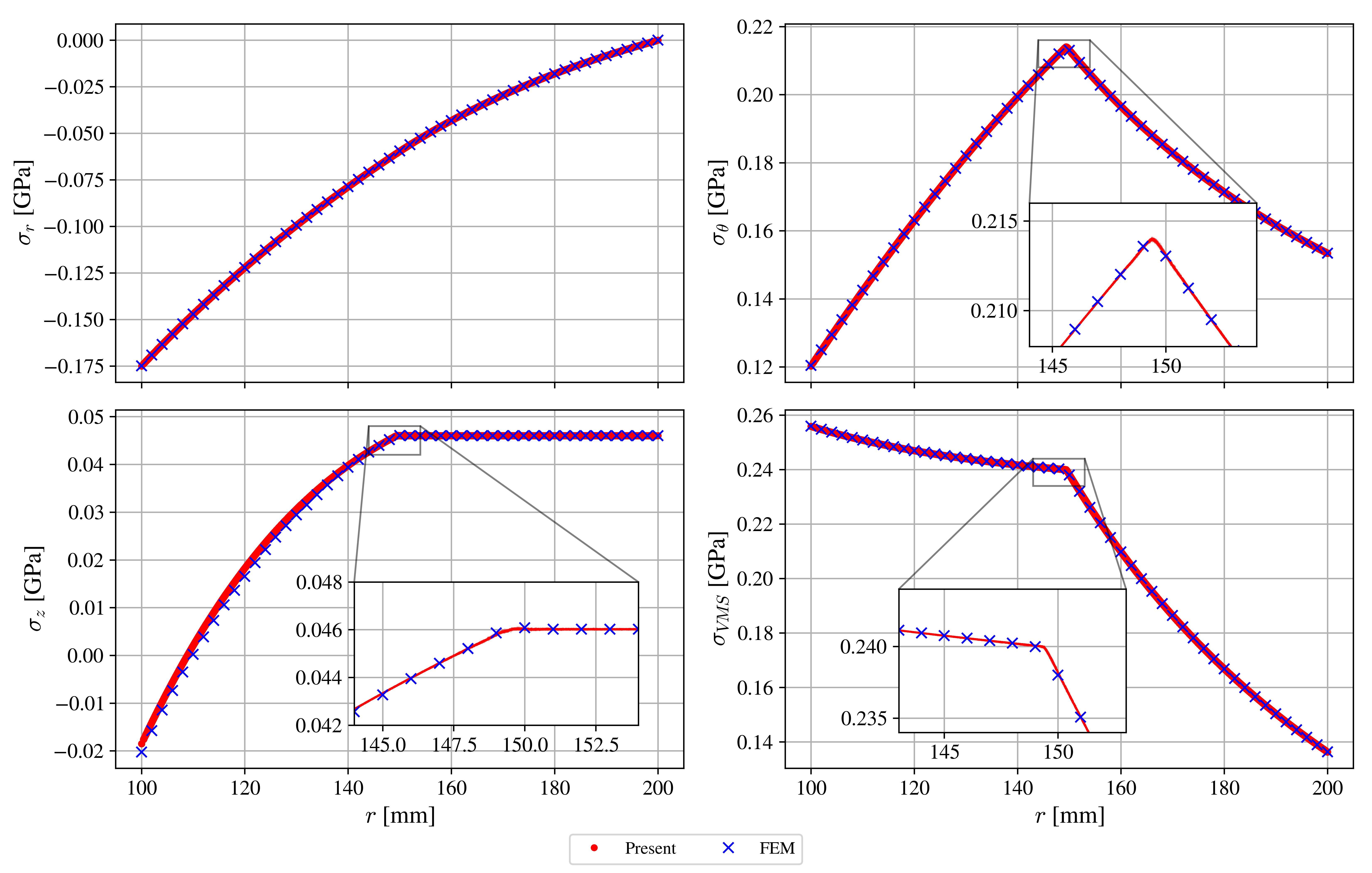}
	\caption{Stress analysis in case of isotropic linear hardening.}
	\label{fig:lenar_hard_stress}
\end{figure}

In the next step, we evaluate the accuracy of the numerical solution in terms of the secondary variables in Figure~\ref{fig:lenar_hard_stress}. Note that, the inside plot of Hoop stress $\sigma_\theta$ in Figure~\ref{fig:lenar_hard_stress} is a close-up to the breaking point, which is positioned at the plastic-elastic interface. All secondary variables are also in good agreement with the FEM-based solution procedure.

A summary of the performance evaluation can be found in Table~\ref{tb: res linear}. We show the von Mises stresses and displacement magnitudes at both extremal values of the domain, i.e.\ $r=a$ and $r=b$. Additionally, the plastic front position $c$ is shown. Note that the same procedure as explained in the previous Section~\ref{sec:perfect_front} has been used to determine the plastic-elastic interface position $c$, with constants $C_i$ gathered in Table~\ref{tab:fit_linear}. 

\begin{table}[H]
	\centering
	\caption{Results of the linear-hardening plastic case - von Mises stress, and displacement magnitude at inner and outer surfaces, and position of the elastic-plastic interface.}
	\label{tb: res linear}
	\begin{threeparttable}[b]
		\def\arraystretch{1.3}
		\begin{tabular}{l|cc|cc|c|c}
			\toprule\toprule
			                                                                                     & \multicolumn{2}{c}{$\svms$ [GPa]} & \multicolumn{2}{|c|}{$\norm{\mathbf u}$ [mm]} & \multirow{2}{*}{$c$ [mm]}   & \multirow{2}{*}{$N$}                                            \\
			                                                                                     & $r = a$                           & $r = b$                                       & $r = a$                     & $r = b$                     &                                   \\
			\hline
			Present                                                                              & 0.25589\tnote{\mjstarsmall}       & 0.13635\tnote{\mjstarsmall}                   & 0.22212\tnote{\mjstarsmall} & 0.13295\tnote{\mjstarsmall} & 149.42\tnote{\mjstar} & $95\,346$ \\
			FEM                                                                                  & 0.25592                           & 0.13634                                       & 0.22224                     & 0.13294                     & 149.40\tnote{\mjstar} & $71\,473$ \\
			\midrule
			$\frac{\widehat{u} - \widehat{u}_\mathrm{FEM}}{\widehat{u}_\mathrm{FEM}} \cdot 10^3$ & -0.117                            & 0.073                                         & -0.540                      & 0.075                       & 0.134                 & -         \\
			\bottomrule
		\end{tabular}
		\begin{tablenotes}
			\item [\mjstarsmall] Average for all nodes on the edge.
			\item [\mjstar] Determined numerically as explaneed in Section~\ref{sec:perfect_front}.
		\end{tablenotes}
	\end{threeparttable}
\end{table}

\begin{table}[H]
	\centering
	\caption{Constant values $C_i$ used to determine the plastic front position $c$ in case of linear hardening.}
	\label{tab:fit_linear}
	\def\arraystretch{1.3}
	\begin{tabular}{l|cccc}
		\toprule\toprule
		regime  & $C_1$                 & $C_2$             & $C_3$                 & $C_4$                \\\midrule
		plastic & $-8.120\cdot 10^{1}$  & $2.426\cdot 10^3$ & $ 1.016\cdot 10^{-1}$ & $1.157$              \\
		elastic & $ 2.993\cdot 10^{-1}$ & $3.055\cdot 10^3$ & $-8.508\cdot 10^{-4}$ & $7.102\cdot 10^{-2}$ \\
		\bottomrule
	\end{tabular}
\end{table}

\section{Additional Example}
\begin{figure}[h]
	\centering
	\includegraphics[width=\textwidth]{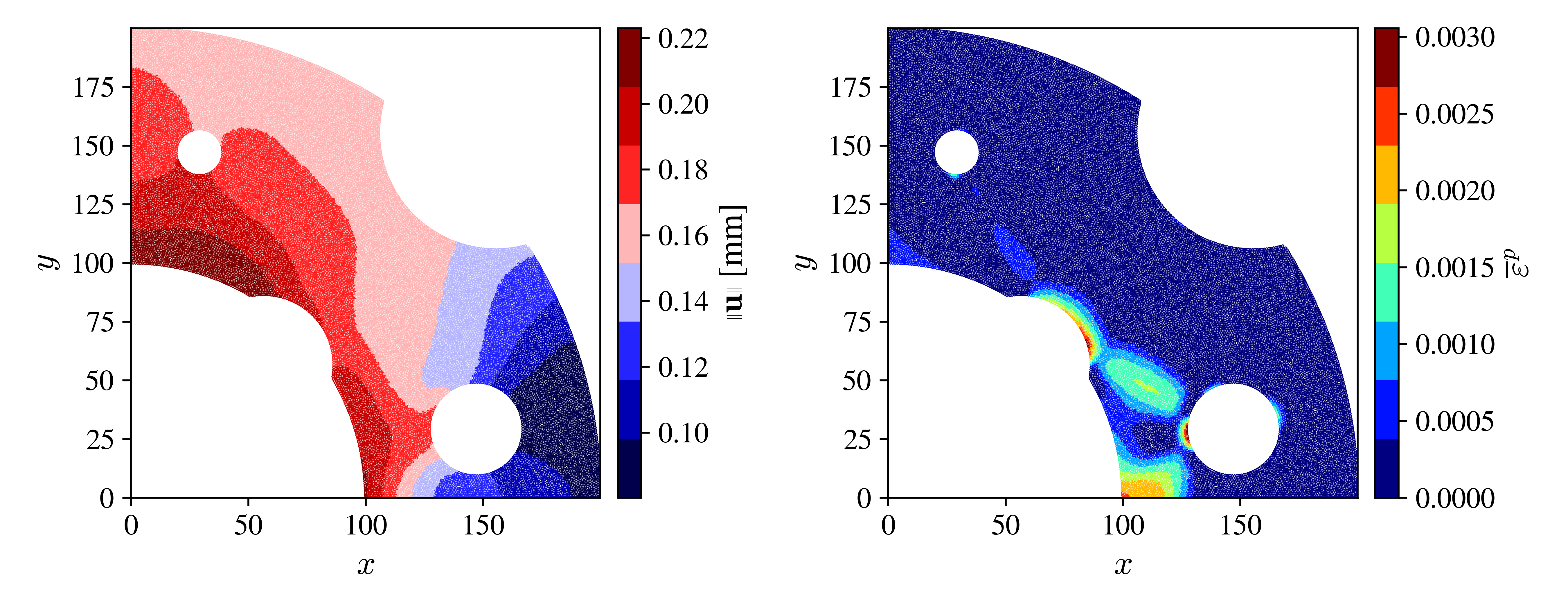}
	\caption{Example solution on irregular domain. Displacement magnitude (left) and accumulated plastic strain (right).}
	\label{fig:irregular}
\end{figure}

In addition to already solved cases, we demonstrate a solution for irregular domain with a more realistic stress-strain relation. Additionally, the domain under consideration is cut with the following four spheres:
\begin{align}
	\Sph 0 & = \left \{\x \in \R^2, \ \left\|\mathbf x - \frac{a+b}{2}\Biggl(\cos\biggl(\frac{\pi}{16}\biggl), \sin\biggl(\frac{\pi}{16}\biggl)\Biggl)\right\| < 20 \right \},   \\
	\Sph 1 & = \left \{\x \in \R^2, \ \left\|\mathbf x - \frac{a+b}{2}\Biggl(\cos\biggl(\frac{7\pi}{16}\biggl), \sin\biggl(\frac{7\pi}{16}\biggl)\Biggl)\right\| < 10 \right \}, \\
	\Sph 2 & = \left \{\x \in \R^2, \ \left\|\mathbf x - 0.8\cdot a \Biggl(\cos\biggl(\frac{\pi}{4}\biggl), \sin\biggl(\frac{\pi}{4}\biggl)\Biggl)\right\| < 30 \right \},       \\
	\Sph 3 & = \left \{\x \in \R^2, \ \left\|\mathbf x - 1.1\cdot b \Biggl(\cos\biggl(\frac{\pi}{4}\biggl), \sin\biggl(\frac{\pi}{4}\biggl)\Biggl)\right\| < 50 \right \},
\end{align}
resulting in a domain as shown in Figure~\ref{fig:irregular}. The edge of $\Sph 2$ is subjected to internal pressure $p=\unit[0.13]{MPa}$, while the remaining three sphere edges are traction free, i.e.\ $\boldsymbol \sigma \vec{n} = 0$.

For the material properties, we simulate a set of uniaxial test results and use a piecewise linear interpolation in the solution procedure. The stress-strain relation is given as provided in Table~\ref{tab:irregular} and illustrated in Figure~\ref{fig:irregular_stress_strain}.

\begin{table}[ht]
	\begin{varwidth}[b]{0.2\linewidth}
	\end{varwidth}
	\hfill
	\begin{varwidth}[b]{0.3\linewidth}
		\centering
		\begin{tabular}{ l | c }
			\toprule
			\toprule
			$\sigma$ [MPa] & $\overline{\varepsilon}^p$ \\
			\midrule
			0.240          & 0.000                      \\
			0.290          & 0.001                      \\
			0.330          & 0.002                      \\
			0.370          & 0.003                      \\
			0.400          & 0.004                      \\
			0.430          & 0.005                      \\
			0.450          & 0.006                      \\
			0.470          & 0.007                      \\
			0.483          & 0.008                      \\
			0.495          & 0.009                      \\
			0.500          & 0.010                      \\
			\bottomrule
		\end{tabular}
		\caption{Synthetic stress-strain relation data.}
		\label{tab:irregular}
	\end{varwidth}
	\hfill
	\begin{varwidth}[b]{0.4\linewidth}
		\centering
		\includegraphics[width=\columnwidth]{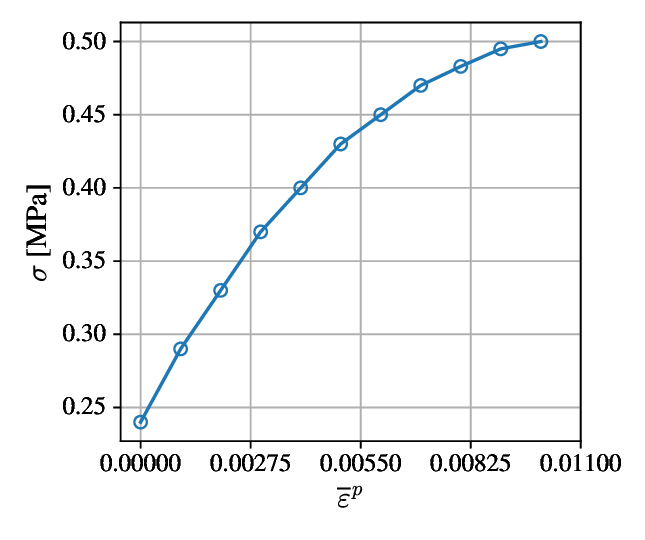}
		\captionof{figure}{Plotted synthetic stress-strain relationship.}
		\label{fig:irregular_stress_strain}
	\end{varwidth}
	\hfill
	\begin{varwidth}[b]{0.1\linewidth}
	\end{varwidth}
\end{table}

The domain is discretized with $20\,292$ scattered nodes corresponding to approximately \unit[1]{mm} internodal spacing. Example solution is shown in Figure~\ref{fig:irregular} with displacement magnitudes on the left and accumulated plastic strain $\overline{\varepsilon}^p$ on the right.
The proposed solution procedure is clearly able to cope with irregular domains and real-like material properties. The largest plastic deformations are observed in the bottom right part of the domain, while the largest deformations are observed in the top right part of the domain, effectively yielding an ellipse-like deformed shape. Note that no special technique or algorithm adjustments were required to solve such problem. The solution procedure used is exactly as described in Algorithm~\ref{alg:adapt}.

\section{Conclusions}
This paper discusses the use of meshless methods, in particular RBF-FD, in solving an elasto-plastic problem using the Picard iteration. The study confirms that the presented solution method effectively accounts for different stages of material behaviour, from elastic deformation to perfect plastic yielding and linear hardening conditions. The method shows good agreement with analytical and FEM solutions (obtained with the commercial Abaqus solver), confirming its accuracy and robustness. It shows stable behaviour and provides results that are independent of the underlying geometry, without the need for special treatment in the approximation of differential operators, boundary conditions or spatial discretisation.

The iterative numerical algorithm initially assumes elastic deformation, i.e. in the first step, displacement and strain fields are calculated by solving the Navier-Cauchy equation. If the stress at any node exceeds the yield stress, indicating plastic deformation, a correction is made using a return mapping algorithm that adjusts the calculated stresses and strains to account for plastic effects, followed by a Picard iteration to solve for non-linear deformations.

Although the Picard iteration converges monotonically and slowly, it avoids the divergence issues common with the Newton-Raphson method, which otherwise converges abruptly, especially when using the consistent tangent matrix~\cite{yarushina2010}. This was also reported in a recent RBF-FD analysis of elasto-plastic benchmark problems~\cite{vuga_improved_2024}, where authors had to introduce additional stabilisations in the calculation of the divergence operator and boundary conditions in order to achieve stable convergence with Newton-Raphson iterations.

In this paper, we show that a direct RBF-FD on pure scattered nodes can be used with the help of Picard iteration, which simplifies both the implementation and the discretisation of the underlying domain considerably, but at the price of slow convergence.

Here, the Navier-Cauchy equation is solved in a strong form using a direct RBF-FD approximation on pure scattered nodes generated by advancing front algorithm~\cite{slak2019generation} that can be generalised to higher dimension~\cite{janvcivc2021monomial} and supports $h$-adaptive discretization~\cite{slak2019rbffd} via spatially variable node spacing.
At each node, the linear differential operators relevant to the problem are approximated using the RBF-FD approximation. Nearby nodes, known as stencil nodes, are used to approximate the differential operators based on RBFs augmented with polynomials. In addition to increased stability and solving problems due to the Gibss phenomenon~\cite{fornberg_solving_2015}, polynomial augmentation also gives us direct control over the order of the method, which can be dynamically adjusted node by node. This property enables a relatively straightforward application of $hp$-adaptivity~\cite{janvcivc2023strong}.

While the local meshless approach offers considerable potential advantages, it is not without drawbacks. Despite extensive research, the RBF-FD method still has some unsolved problems. Node quality in a meshless context is not nearly as well understood as in FEM, where mesh quality measures such as aspect ratio, skewness, Jacobian ratio, etc. are well-defined concepts~\cite{gokhale2008practical}. There are also no clear guidelines for the selection of stencils. In this work, we simply used the nearest $n$ nodes as a stencil, which leads to a relatively computationally intensive method due to the rather high stencil size reuqirements~\cite{bayona2019insight}.

Future work could focus on extending this method to three-dimensional problems and investigating the effects of dynamic loading conditions in more complex hardening models. As for the approximation method, symmetric stencils~\cite{davydov_improved_2023} could be used in addition to simple closest neighbourhood stencils, to improve the execution performance. The presented method is also compatible with the recently introduced $h$~\cite{slak2019rbffd} and $hp$~\cite{janvcivc2023strong} adaptivities, which could be used to improve the performance of elasto-plastic deformation simulation. As for the iteration, an improved Picard iteration could be considered, which has been recently applied to flow in porous media~\cite{ZHU2024105979}. The presented solution could also be coupled with a more sophisticated Newton-Raphson iteration to alleviate its stability problems~\cite{yarushina2010}. 

\section*{Acknowledgments}
Authors acknowledge the financial support from the Slovenian Research And Innovation Agency (ARIS) research core funding No.\ P2-0095.

\bibliography{refs}

\begin{thebibliography}{10}
\expandafter\ifx\csname url\endcsname\relax
  \def\url#1{\texttt{#1}}\fi
\expandafter\ifx\csname urlprefix\endcsname\relax\def\urlprefix{URL }\fi
\expandafter\ifx\csname href\endcsname\relax
  \def\href#1#2{#2} \def\path#1{#1}\fi

\bibitem{chakrabarty1987}
J.~Chakrabarty, Theory of Plasticity, McGraw-Hill international editions.
  Engineering mechanics series, McGraw-Hill, 1987.

\bibitem{hill1998}
R.~Hill, The Mathematical Theory of Plasticity, Oxford classic texts in the
  physical sciences, Clarendon Press, 1998.

\bibitem{desouza2008}
E.~A. de~Souza~Neto, D.~Peric, D.~R.~J. Owen, Computational Methods for
  Plasticity: Theory and Applications, Wiley, 2008.

\bibitem{Schrder2015SmallSP}
B.~Schr\"{o}der, D.~Kuhl, Small strain plasticity: classical versus multifield
  formulation, Archive of Applied Mechanics 85 (2015) 1127--1145.

\bibitem{Roostaei2018ACS}
A.~A. Roostaei, H.~Jahed, A cyclic small-strain plasticity model for wrought mg
  alloys under multiaxial loading: Numerical implementation and validation,
  International Journal of Mechanical Sciences.

\bibitem{amouzou2021}
G.~Y. Amouzou, A.~Soula\"{i}mani, Numerical algorithms for elastoplacity:
  Finite elements code development and implementation of the mohr–coulomb
  law, Applied Sciences 11~(10).
\newblock \href {http://dx.doi.org/10.3390/app11104637}
  {\path{doi:10.3390/app11104637}}.

\bibitem{mitchell2014comparison}
W.~F. Mitchell, M.~A. McClain, A comparison of hp-adaptive strategies for
  elliptic partial differential equations, ACM Transactions on Mathematical
  Software (TOMS) 41~(1) (2014) 1--39.

\bibitem{segeth2010review}
K.~Segeth, A review of some a posteriori error estimates for adaptive finite
  element methods, Mathematics and Computers in Simulation 80~(8) (2010)
  1589--1600.

\bibitem{cottrell2009isogeometric}
J.~A. Cottrell, T.~J. Hughes, Y.~Bazilevs, Isogeometric analysis: toward
  integration of CAD and FEA, John Wiley \& Sons, 2009.

\bibitem{liu_introduction_2005}
G.-R. Liu, Y.-T. Gu, An introduction to meshfree methods and their programming,
  Springer Science \& Business Media, 2005.

\bibitem{liu2002mesh}
G.-R. Liu, Mesh free methods: moving beyond the finite element method, CRC
  press, 2002.
\newblock \href {http://dx.doi.org/10.1201/9781420040586}
  {\path{doi:10.1201/9781420040586}}.

\bibitem{Reuther2012}
K.~Reuther, B.~Sarler, M.~Rettenmayr, Solving diffusion problems on an
  unstructured, amorphous grid by a meshless method, Int. J. Therm. Sci. 51
  (2012) 16--22.
\newblock \href {http://dx.doi.org/10.1016/j.ijthermalsci.2011.08.017}
  {\path{doi:10.1016/j.ijthermalsci.2011.08.017}}.

\bibitem{wendland2004scattered}
H.~Wendland, Scattered data approximation, Vol.~17, Cambridge university press,
  2004.

\bibitem{liu2009meshfree}
G.-R. Liu, Meshfree methods: moving beyond the finite element method, CRC
  press, 2009.

\bibitem{slak2019generation}
J.~Slak, G.~Kosec, On generation of node distributions for meshless pde
  discretizations, SIAM Journal on Scientific Computing 41~(5) (2019)
  A3202--A3229.
\newblock \href {http://dx.doi.org/10.1137/18M1231456}
  {\path{doi:10.1137/18M1231456}}.

\bibitem{duh2021fast}
U.~Duh, G.~Kosec, J.~Slak, Fast variable density node generation on parametric
  surfaces with application to mesh-free methods, SIAM Journal on Scientific
  Computing 43~(2) (2021) A980--A1000.

\bibitem{fornberg2015fast}
B.~Fornberg, N.~Flyer, Fast generation of 2-{D} node distributions for
  mesh-free {PDE} discretizations, Computers \& Mathematics with Applications
  69~(7) (2015) 531–544.
\newblock \href {http://dx.doi.org/10.1016/j.camwa.2015.01.009}
  {\path{doi:10.1016/j.camwa.2015.01.009}}.

\bibitem{van2021fast}
K.~van~der Sande, B.~Fornberg, Fast variable density 3-d node generation, SIAM
  Journal on Scientific Computing 43~(1) (2021) A242--A257.

\bibitem{shankar2018robust}
V.~Shankar, R.~M. Kirby, A.~L. Fogelson, Robust node generation for meshfree
  discretizations on irregular domains and surfaces, SIAM J. Sci. Comput.
  40~(4) (2018) 2584--2608.
\newblock \href {http://dx.doi.org/10.1137/17m114090x}
  {\path{doi:10.1137/17m114090x}}.

\bibitem{duh2023discretization}
U.~Duh, V.~Shankar, G.~Kosec, Discretization of non-uniform rational b-spline
  (nurbs) models for meshless isogeometric analysis, arXiv preprint
  arXiv:2303.02638.

\bibitem{ji2005meshfree}
Z.~Ji-fa, Z.~Wen-pu, Z.~Yao, A meshfree method and its applications to
  elasto-plastic problems, Journal of Zhejiang University-SCIENCE A 6~(2)
  (2005) 148--154.

\bibitem{kargarnovin2004elasto}
M.~Kargarnovin, H.~E. Toussi, S.~Fariborz, Elasto-plastic element-free galerkin
  method, Computational Mechanics 33 (2004) 206--214.

\bibitem{app132312591}
J.~Belinha, M.~Aires, D.~E. Rodrigues,
  \href{https://www.mdpi.com/2076-3417/13/23/12591}{Elastoplastic analysis of
  frame structures using radial point interpolation meshless methods}, Applied
  Sciences 13~(23).
\newblock \href {http://dx.doi.org/10.3390/app132312591}
  {\path{doi:10.3390/app132312591}}.
\newline\urlprefix\url{https://www.mdpi.com/2076-3417/13/23/12591}

\bibitem{kosec_weak_2019}
G.~Kosec, J.~Slak, M.~Depolli, R.~Trobec, K.~Pereira, S.~Tomar, T.~Jacquemin,
  S.~P. Bordas, M.~A. Wahab, Weak and strong from meshless methods for linear
  elastic problem under fretting contact conditions, Tribology
  InternationalPublisher: Elsevier.
\newblock \href {http://dx.doi.org/10.1016/j.triboint.2019.05.041}
  {\path{doi:10.1016/j.triboint.2019.05.041}}.

\bibitem{JANKOWSKA201812}
M.~A. Jankowska, On elastoplastic analysis of some plane stress problems with
  meshless methods and successive approximations method, Engineering Analysis
  with Boundary Elements 95 (2018) 12--24.
\newblock \href
  {http://dx.doi.org/https://doi.org/10.1016/j.enganabound.2018.06.013}
  {\path{doi:https://doi.org/10.1016/j.enganabound.2018.06.013}}.

\bibitem{jiang2021nonlinear}
P.~Jiang, H.~Zheng, J.~Xiong, P.~Wen, Nonlinear elastic-plastic analysis of
  reinforced concrete column-steel beam connection by rbf-fd method,
  Engineering Analysis with Boundary Elements 128 (2021) 188--194.

\bibitem{bayona2017role}
V.~Bayona, N.~Flyer, B.~Fornberg, G.~A. Barnett, On the role of polynomials in
  rbf-fd approximations: Ii. numerical solution of elliptic pdes, Journal of
  Computational Physics 332 (2017) 257--273.

\bibitem{wang2002optimal}
J.~Wang, G.~Liu, On the optimal shape parameters of radial basis functions used
  for 2-d meshless methods, Computer methods in applied mechanics and
  engineering 191~(23-24) (2002) 2611--2630.

\bibitem{vuga_improved_2024}
G.~Vuga, B.~Mavrič, B.~Šarler,
  \href{https://linkinghub.elsevier.com/retrieve/pii/S0045782523006254}{An
  improved local radial basis function method for solving small-strain
  elasto-plasticity}, Computer Methods in Applied Mechanics and Engineering 418
  (2024) 116501.
\newblock \href {http://dx.doi.org/10.1016/j.cma.2023.116501}
  {\path{doi:10.1016/j.cma.2023.116501}}.
\newline\urlprefix\url{https://linkinghub.elsevier.com/retrieve/pii/S0045782523006254}

\bibitem{yarushina2010}
V.~M. Yarushina, M.~Dabrowski, Y.~Y. Podladchikov, An analytical benchmark with
  combined pressure and shear loading for elastoplastic numerical models,
  Geochemistry, Geophysics, Geosystems 11~(8).
\newblock \href {http://dx.doi.org/https://doi.org/10.1029/2010GC003130}
  {\path{doi:https://doi.org/10.1029/2010GC003130}}.

\bibitem{KOLODZIEJ20134217}
J.~A. Kolodziej, M.~A. Jankowska, M.~Mierzwiczak,
  \href{https://www.sciencedirect.com/science/article/pii/S0020768313003417}{Meshless
  methods for the inverse problem related to the determination of elastoplastic
  properties from the torsional experiment}, International Journal of Solids
  and Structures 50~(25) (2013) 4217--4225.
\newblock \href
  {http://dx.doi.org/https://doi.org/10.1016/j.ijsolstr.2013.08.025}
  {\path{doi:https://doi.org/10.1016/j.ijsolstr.2013.08.025}}.
\newline\urlprefix\url{https://www.sciencedirect.com/science/article/pii/S0020768313003417}

\bibitem{moayyedian2021elastic}
F.~Moayyedian, J.~K. Grabski, Elastic--plastic torsion problem with non-linear
  hardenings using the method of fundamental solution, Archives of Civil and
  Mechanical Engineering 21~(4) (2021) 155.

\bibitem{XU2023939}
B.~Xu, R.~Zhang, K.~Yang, G.~Yu, Y.~Chen,
  \href{https://www.sciencedirect.com/science/article/pii/S0955799722004295}{Application
  of generalized finite difference method for elastoplastic torsion analysis of
  prismatic bars}, Engineering Analysis with Boundary Elements 146 (2023)
  939--950.
\newblock \href
  {http://dx.doi.org/https://doi.org/10.1016/j.enganabound.2022.11.028}
  {\path{doi:https://doi.org/10.1016/j.enganabound.2022.11.028}}.
\newline\urlprefix\url{https://www.sciencedirect.com/science/article/pii/S0955799722004295}

\bibitem{slak2021medusa}
J.~Slak, G.~Kosec, Medusa: A c++ library for solving pdes using strong form
  mesh-free methods, ACM Transactions on Mathematical Software (TOMS) 47~(3)
  (2021) 1--25.

\bibitem{janvcivc2021monomial}
M.~Jan{\v{c}}i{\v{c}}, J.~Slak, G.~Kosec, Monomial augmentation guidelines for
  rbf-fd from accuracy versus computational time perspective, Journal of
  Scientific Computing 87~(1) (2021) 9.

\bibitem{janvcivc2023strong}
M.~Jan{\v{c}}i{\v{c}}, G.~Kosec, Strong form mesh-free hp-adaptive solution of
  linear elasticity problem, Engineering with Computers (2023) 1--21.

\bibitem{najafi_divergence-free_2022}
M.~Najafi, M.~Dehghan, B.~Šarler, G.~Kosec, B.~Mavrič,
  \href{https://link.springer.com/10.1007/s00366-022-01621-w}{Divergence-free
  meshless local {Petrov}–{Galerkin} method for {Stokes} flow}, Engineering
  with Computers 38~(6) (2022) 5359--5377.
\newblock \href {http://dx.doi.org/10.1007/s00366-022-01621-w}
  {\path{doi:10.1007/s00366-022-01621-w}}.
\newline\urlprefix\url{https://link.springer.com/10.1007/s00366-022-01621-w}

\bibitem{berljavac_rbf-fd_2021}
J.~M. Berljavac, P.~K. Mishra, J.~Slak, G.~Kosec, {RBF}-{FD} analysis of {2D}
  time-domain acoustic wave propagation in heterogeneous media, Computers \&
  Geosciences 153 (2021) 104796, publisher: Elsevier.

\bibitem{jancic_meshless_2024}
M.~Jančič, M.~Založnik, G.~Kosec,
  \href{https://linkinghub.elsevier.com/retrieve/pii/S0021999124002225}{Meshless
  interface tracking for the simulation of dendrite envelope growth}, Journal
  of Computational Physics 507 (2024) 112973.
\newblock \href {http://dx.doi.org/10.1016/j.jcp.2024.112973}
  {\path{doi:10.1016/j.jcp.2024.112973}}.
\newline\urlprefix\url{https://linkinghub.elsevier.com/retrieve/pii/S0021999124002225}

\bibitem{slak2019rbffd}
J.~Slak, G.~Kosec, Adaptive radial basis function-generated finite differences
  method for contact problems, International Journal for Numerical Methods in
  Engineering 119~(7) (2019) 661--686.
\newblock \href {http://dx.doi.org/https://doi.org/10.1002/nme.6067}
  {\path{doi:https://doi.org/10.1002/nme.6067}}.

\bibitem{depolli_parallel_2022}
M.~Depolli, J.~Slak, G.~Kosec,
  \href{https://linkinghub.elsevier.com/retrieve/pii/S0045794922000335}{Parallel
  domain discretization algorithm for {RBF}-{FD} and other meshless numerical
  methods for solving {PDEs}}, Computers \& Structures 264 (2022) 106773.
\newline\urlprefix\url{https://linkinghub.elsevier.com/retrieve/pii/S0045794922000335}

\bibitem{davydov2011adaptive}
O.~Davydov, D.~T. Oanh, Adaptive meshless centres and rbf stencils for poisson
  equation, Journal of Computational Physics 230~(2) (2011) 287--304.

\bibitem{davydov2023improved}
O.~Davydov, D.~T. Oanh, N.~M. Tuong, Improved stencil selection for meshless
  finite difference methods in 3d, Journal of Computational and Applied
  Mathematics (2023) 115031.

\bibitem{bayona2019comparison}
V.~Bayona, Comparison of moving least squares and rbf+ poly for interpolation
  and derivative approximation, Journal of Scientific Computing 81 (2019)
  486--512.

\bibitem{bayona2019insight}
V.~Bayona, An insight into rbf-fd approximations augmented with polynomials,
  Computers \& Mathematics with Applications 77~(9) (2019) 2337--2353.

\bibitem{flyer2016role}
N.~Flyer, B.~Fornberg, V.~Bayona, G.~A. Barnett, On the role of polynomials in
  rbf-fd approximations: I. interpolation and accuracy, Journal of
  Computational Physics 321 (2016) 21--38.

\bibitem{flyer2016}
N.~Flyer, B.~FornXberg, V.~Bayona, G.~A. Barnett, On the role of polynomials in
  rbf-fd approximations: I. interpolation and accuracy, Journal of
  Computational Physics 321 (2016) 21--38.
\newblock \href {http://dx.doi.org/https://doi.org/10.1016/j.jcp.2016.05.026}
  {\path{doi:https://doi.org/10.1016/j.jcp.2016.05.026}}.

\bibitem{le2023guidelines}
S.~Le~Borne, W.~Leinen, Guidelines for rbf-fd discretization: Numerical
  experiments on the interplay of a multitude of parameter choices, Journal of
  Scientific Computing 95~(1) (2023) 8.

\bibitem{9803334}
M.~Jančič, G.~Kosec, Stability analysis of rbf-fd and wls based local strong
  form meshless methods on scattered nodes, in: 2022 45th Jubilee International
  Convention on Information, Communication and Electronic Technology (MIPRO),
  2022, pp. 275--280.
\newblock \href {http://dx.doi.org/10.23919/MIPRO55190.2022.9803334}
  {\path{doi:10.23919/MIPRO55190.2022.9803334}}.

\bibitem{SLAK20193}
J.~Slak, G.~Kosec, Refined meshless local strong form solution of cauchy-navier
  equation on an irregular domain, Engineering Analysis with Boundary Elements
  100 (2019) 3--13, improved Localized and Hybrid Meshless Methods - Part 1.
\newblock \href
  {http://dx.doi.org/https://doi.org/10.1016/j.enganabound.2018.01.001}
  {\path{doi:https://doi.org/10.1016/j.enganabound.2018.01.001}}.

\bibitem{eigenweb}
G.~Guennebaud, B.~Jacob, et~al., Eigen v3, http://eigen.tuxfamily.org (2010).

\bibitem{krabbenhoft2002basic}
K.~Krabbenh{\o}ft, Basic computational plasticity, University of Denmark.

\bibitem{PEREZFOGUET20014627}
A.~Perez-Foguet, A.~Rodriguez-Ferran, A.~Huerta,
  \href{https://www.sciencedirect.com/science/article/pii/S0045782500003364}{Consistent
  tangent matrices for substepping schemes}, Computer Methods in Applied
  Mechanics and Engineering 190~(35) (2001) 4627--4647.
\newblock \href
  {http://dx.doi.org/https://doi.org/10.1016/S0045-7825(00)00336-4}
  {\path{doi:https://doi.org/10.1016/S0045-7825(00)00336-4}}.
\newline\urlprefix\url{https://www.sciencedirect.com/science/article/pii/S0045782500003364}

\bibitem{ROOSTAEI202223}
A.~A. Roostaei, H.~Jahed, 2 - fundamentals of cyclic plasticity models, in:
  H.~Jahed, A.~A. Roostaei (Eds.), Cyclic Plasticity of Metals, Elsevier Series
  on Plasticity of Materials, Elsevier, 2022, pp. 23--51.
\newblock \href
  {http://dx.doi.org/https://doi.org/10.1016/B978-0-12-819293-1.00011-5}
  {\path{doi:https://doi.org/10.1016/B978-0-12-819293-1.00011-5}}.

\bibitem{fornberg_solving_2015}
B.~Fornberg, N.~Flyer,
  \href{https://www.cambridge.org/core/product/identifier/S0962492914000130/type/journal_article}{Solving
  {PDEs} with radial basis functions}, Acta Numerica 24 (2015) 215--258.
\newblock \href {http://dx.doi.org/10.1017/S0962492914000130}
  {\path{doi:10.1017/S0962492914000130}}.
\newline\urlprefix\url{https://www.cambridge.org/core/product/identifier/S0962492914000130/type/journal_article}

\bibitem{gokhale2008practical}
N.~S. Gokhale, Practical finite element analysis, Finite to infinite, 2008.

\bibitem{davydov_improved_2023}
O.~Davydov, D.~T. Oanh, N.~M. Tuong,
  \href{https://linkinghub.elsevier.com/retrieve/pii/S037704272200629X}{Improved
  stencil selection for meshless finite difference methods in {3D}}, Journal of
  Computational and Applied Mathematics 425 (2023) 115031.
\newblock \href {http://dx.doi.org/10.1016/j.cam.2022.115031}
  {\path{doi:10.1016/j.cam.2022.115031}}.
\newline\urlprefix\url{https://linkinghub.elsevier.com/retrieve/pii/S037704272200629X}

\bibitem{ZHU2024105979}
S.~Zhu, L.~Zhang, L.~Wu, L.~Tan, H.~Chen,
  \href{https://www.sciencedirect.com/science/article/pii/S0266352X2300736X}{Application
  of improved picard iteration method to simulate unsaturated flow and
  deformation in deformable porous media}, Computers and Geotechnics 166 (2024)
  105979.
\newblock \href
  {http://dx.doi.org/https://doi.org/10.1016/j.compgeo.2023.105979}
  {\path{doi:https://doi.org/10.1016/j.compgeo.2023.105979}}.
\newline\urlprefix\url{https://www.sciencedirect.com/science/article/pii/S0266352X2300736X}

\end{thebibliography}
\end{document}